\documentclass[11pt]{amsart}
\usepackage{amssymb}
\usepackage{amsthm}

\newtheorem{theorem}{Theorem}[section]

\newtheorem{corollary}[theorem]{Corollary}
\newtheorem{lemma}[theorem]{Lemma}
\newtheorem{proposition}[theorem]{Proposition}
\theoremstyle{definition}
\newtheorem{definition}[theorem]{Definition}

\newtheorem{question}[theorem]{Question}
\theoremstyle{remark}
\newtheorem{remark}[theorem]{Remark}
\newtheorem{example}[theorem]{Example}

\newtheorem{claim}{\bf Claim}
\numberwithin{equation}{section}

\def\R {{\mathbb{R}}}
\def\T {{\mathbb{T}}}

\def\Z {{\mathbb{Z}}}

\def\3{{|\!|\!|}}

\begin{document}

\title[On equivalence relations induced by TSI Polish groups]{On equivalence relations induced by Polish groups admitting compatible two-sided invariant metrics}
\author{Longyun Ding}
\address{School of Mathematical Sciences and LPMC, Nankai University, Tianjin, 300071, P.R.China}
\email{dingly@nankai.edu.cn}
\thanks{This research is partially supported by the National Natural Science Foundation of China (Grant No. 12271264).}
\author{Yang Zheng}
\email{1120200015@mail.nankai.edu.cn}

\subjclass[2010]{03E15, 22A05, 22D05, 22E15, 46A16}
\keywords{Borel reduction, TSI Polish group, equivalence relation}


\begin{abstract}
Given a Polish group $G$, let $E(G)$ be the right coset equivalence relation $G^\omega/c(G)$, where $c(G)$ is the group of all convergent sequences in $G$. We first established two results:
\begin{enumerate}
\item Let $G,H$ be two Polish groups. If $H$ is TSI but $G$ is not, then $E(G)\not\le_BE(H)$.
\item Let $G$ be a Polish group. Then the following are equivalent: (a) $G$ is TSI non-archimedean; (b)$E(G)\leq_B E_0^\omega$; and (c) $E(G)\leq_B \R^\omega/c_0$.
In particular, $E(G)\sim_B E_0^\omega$ iff $G$ is TSI uncountable non-archimedean.
\end{enumerate}

A critical theorem presented in this article is as follows: Let $G$ be a TSI Polish group, and let $H$ be a closed subgroup of the product of a sequence of TSI strongly NSS Polish groups. If $E(G)\le_BE(H)$, then there exists a continuous homomorphism $S:G_0\rightarrow H$ such that $\ker(S)$ is non-archimedean, where $G_0$ is the connected component of the identity of $G$. The converse holds if $G$ is connected, $S(G)$ is closed in $H$, and the interval $[0,1]$ can be embedded into $H$.

As its applications, we prove several Rigid theorems for TSI Lie groups, locally compact Polish groups, separable Banach spaces, and separable Fr\'echet spaces, respectively.
\end{abstract}
\maketitle

\section{Introduction}

In recent years, logicians have achieved remarkable research outcomes in descriptive set theory, specifically in the investigation of the relative complexity among equivalence relations originating from various branches of mathematics, utilizing Borel reducibility. Polish groups and their actions play a crucial role in this research direction. Concurrently, researchers also aim to employ Borel reducibility among equivalence relations to characterize the properties of Polish groups.

The authors introduced in~\cite{DZ} the notion of \textit{equivalence relations induced by Polish groups}: given a Polish group $G$, the equivalence relation $E(G)$ is defined on $G^\omega$ as:
$$ xE(G)y \iff \lim_n x(n)y(n)^{-1} \mbox{ converges in } G,$$
for $x,y\in G^\omega$.
These equivalence relations have shown great potential in characterizing properties of Polish groups. In fact, based on the results on Borel reducibility involving the equivalence relation $E(G)$, we can accurately determine the classes to which certain Polish groups $G$ belong.
For instance, (1) $G$ is countable discrete iff $E(G)\sim_B E_0$; (2) $G$ is non-archimedean iff $E(G)\le_B=^+$; and (3) if $H$ is CLI but $G$ is not, then $E(G)\nleq_BE(H)$ (see~\cite{DZ}).
Additionally, the authors further established results in~\cite{DZlcoabe} on Borel reducibility among equivalence relations induced by locally compact abelian Polish groups, including a Rigid Theorem:

\begin{theorem}[Rigid Theorem,~{\cite[Theorem 2.8]{DZlcoabe}}]\label{abelian}
Let $G$ be a compact connected abelian Polish group and $H$ a locally compact abelian Polish group. Then $E(G)\le_BE(H)$ iff there is a continuous homomorphism $S:G\to H$ such that $\ker(S)$ is non-archimedean.
\end{theorem}

In this article, we shift our focus to \textit{TSI Polish groups}, i.e., those that admit compatible complete two-sided invariant metrics. Firstly, the Borel reducibility among equivalence relations induced by Polish groups can accurately distinguish non-TSI and TSI Polish groups. Actually, we have concluded that:

\begin{theorem}
Let $G,H$ be two Polish groups. If $H$ is TSI but $G$ is not, then $E(G)\not\le_BE(H)$.
\end{theorem}

To elucidate the research significance of equivalence relations induced by Polish groups, we compare them with benchmark equivalence relations such as $E_0$, $E_0^\omega$, $\R^\omega/c_0$, etc. (definitions of these benchmark equivalence relations can be found in the next section).

Surprisingly, we prove that there is NO Polish group $G$ such that
$$E_0^\omega<_B E(G)\leq_B \R^\omega/ c_0.$$
In stark contrast, Farah proved that the partially ordered set $P(\omega)/{\rm Fin}$ can be embedded into Borel equivalence relations between $E_0^\omega$ and $\R^\omega/c_0$ (see~\cite[Theorem 5.4]{farah}).

The above results stem from the following more precise theorem.

\begin{theorem}
Let $G$ be a Polish group. Then the following are equivalent:
\begin{enumerate}
  \item $G$ is TSI non-archimedean;
  \item $E(G)\leq_B E_0^\omega$; and
  \item $E(G)\leq_B \R^\omega/c_0$.
\end{enumerate}
In particular, $E(G)\sim_B E_0^\omega$ iff $G$ is TSI uncountable non-archimedean.
\end{theorem}

To generalize the Rigid Theorem for locally compact abelian Polish groups mentioned above, we shall first generalize \cite[Theorem 6.13]{DZ} to a highly technical theorem, namely the Pre-rigid Theorem (the statement of this theorem is too lengthy to include in the introduction).
The Pre-rigid Theorem and all its applications involve a notion named strongly NSS Polish groups.
A Polish group $G$ is called \textit{strongly NSS} if there exists an open neighborhood $V$ of $1_G$ in $G$ such that
$$\forall(g_n)\in G^\omega\,(\lim_n g_n\ne 1_G\Rightarrow\exists n_0<\cdots<n_k\,(g_{n_0}\cdots g_{n_k}\notin V)).$$

We employ the Pre-rigid Theorem to prove the following theorem, where $G_0$ is the connected component of $1_G$ in $G$.

\begin{theorem}
Let $G$ be a TSI Polish group, and let $H$ be a closed subgroup of the product of a sequence of TSI strongly NSS Polish groups. If $E(G)\le_BE(H)$, then there exists a continuous homomorphism $S:W\rightarrow H$ such that $\ker(S)$ is non-archimedean, where $W\supseteq G_0$ is a countable intersection of clopen subgroups in $G$.

In particular, the converse holds if $G=W$, $S(G)$ is closed in $H$, and the interval $[0,1]$ can be embedded into $H$.
\end{theorem}

Several applications of the Pre-rigid Theorem and the above theorem are listed below.

\begin{theorem}
Let $G$ be a TSI Polish group, and let $H$ be a closed subgroup of the product of a sequence of TSI strongly NSS Polish groups. If $H$ is totally disconnected but $G$ is not, then $E(G)\nleq_BE(H)$.
\end{theorem}

Recall that a \textit{Lie group} is a group which is also a smooth manifold such that the group operations are smooth functions.
Let $G$ be a Lie group, then $G_0$ is an open normal subgroup of $G$. A completely metrizable topological group $G$ is called a \textit{pro-Lie group} if every open neighborhood of $1_G$ contains a normal subgroup $N$ such that $G/N$ is a Lie group~(cf.~\cite[Definition 1]{HM07}).

\begin{theorem}
Let $G,H$ be two TSI Polish groups such that $H$ is a pro-Lie group. If $E(G)\leq_B E(H)$, then there exists a continuous homomorphism $S:G_0\rightarrow H$ such that $\ker(S)$ is non-archimedean.

In particular, the converse holds if $G$ is connected and $S(G)$ is closed in $H$.
\end{theorem}

The following theorem is a generalization of Theorem~\ref{abelian}. It should be emphasized that all locally compact TSI groups are pro-Lie groups.

\begin{theorem}[Rigid Theorem for locally compact TSI groups]
Let $G$ be a locally compact connected TSI Polish group and $H$ a TSI pro-Lie Polish group. Then $E(G)\leq_B E(H)$ iff there exists a continuous homomorphism $S:G\rightarrow H$ such that $\ker(S)$ is non-archimedean.
\end{theorem}

Note that a topological group $G$ is a Lie group iff it is locally compact and there exists an open neighborhood $V$ of $1_G$ in $G$ such that no non-trivial subgroup of $G$ is contained in $V$. This leads to a positive answer to Question 7.4 of~\cite{DZ} as follows:

\begin{theorem}[Rigid Theorem for TSI Lie groups]
Let $G,H$ be two separable TSI Lie groups such that $G$ is connected. Then $E(G)\leq_B E(H)$ iff there exists a continuous locally injective homomorphism $S:G\rightarrow H$.
\end{theorem}

All separable \textit{Fr{\'e}chet spaces}, i.e., separable completely  metrizable topological vector spaces, can be viewed as abelian Polish groups under the addition operation.

\begin{theorem}[Rigid Theorem for Fr\'echet spaces]
Let $X,Y$ be two separable Fr{\'e}chet spaces such that $Y$ is a closed subgroup of the product of a sequence of TSI strongly NSS Polish groups. Then $E(X)\leq_B E(Y)$ iff $X$ is topologically isomorphic to a closed linear subspace of $Y$.
\end{theorem}

All Banach spaces are Fr\'echet spaces. We point out that a separable Banach space $X$ is not strongly NSS iff it has a closed linear subspace topologically isomorphic to $c_0$.
This implies that:

\begin{theorem}[Rigid Theorem for Banach spaces]
Let $X,Y$ be two separable Banach spaces such that $Y$ contains no closed linear subspaces topologically isomorphic to $c_0$. Then $E(X)\leq_B E(Y)$ iff $X$ is topologically isomorphic to a closed linear subspace of $Y$.
\end{theorem}

In addition, we attempt to study some examples induced by totally disconnected TSI Polish groups.
For $p\in[1,+\infty)$ and $a\in c_0$, let
$$I_a=\{n\in\omega: a(n)\neq 0\},\quad A_{p,a}=\{v\in l_p:\forall n\,(v(n)\in a(n)\mathbb{Z})\}.$$
If $I_a$ is infinite, then $A_{p,a}$ is a totally disconnected, strongly NSS abelian Polish group, but is not non-archimedean.
In particular, we let $d(n)=2^{-n}$, and define
$$A_p=A_{p,d}=\{v\in l_p:\forall n\,(v(n)\in 2^{-n}\mathbb Z)\}.$$

\begin{theorem}
For $p,q\in[1,+\infty)$ and $a\in c_0$, the following hold:
\begin{enumerate}
  \item [(1)] if $I_a$ is a nonempty finite set, then $E(A_{p,a})\sim_B E_0$;
  \item [(2)] if $I_a$ is infinite, then $E(A_{p,a})\sim_B E(A_p)$;
  \item [(3)] $E(A_p)<_B E(l_q)$;
  \item [(4)] $E(A_p)\leq_B E(l_q)\iff p=q$; and
  \item [(5)] $E(A_p)\leq_B E(A_q)\iff p=q$.
\end{enumerate}
\end{theorem}

This article is organized as follows. In section 2, we recall some notions in descriptive set theory and also recall some notions and results originated from~\cite{DZ} that will be repeatedly used in this article.
In section 3, we prove Theorem~1.2. It is worth noting that some notation defined in this section will continue to be used in Section 5.
In section 4, we prove Theorem~1.3.
In section 5, we prove the Pre-rigid Theorem, which will serve as the foundation for the subsequent sections.
In section 6, we prove theorems~1.4--1.10.
Finally, in section 7, we prove Theorem~1.11.

\section{Preliminaries}

In this article, all groups are assumed to contain at least two elements.
Any linear space can be viewed as an abelian group. The addition operation in it, as well as in all its subgroups, is denoted by $+$, and its identity element is denoted by $0$.
Unless otherwise specified, for any abstract topological group $G$, we use multiplicative notation to express the group operation, and $1_G$ to express the identity element of $G$.

Given a topological group $G$, the connected component of $1_G$ is denoted by $G_0$, which is clearly a closed normal subgroup of $G$.
Note that $G$ is \textit{totally disconnected} iff $G_0=\{1_G\}$.
It is worth noting that any open subgroup $H$ of $G$ is also closed, since $H=G\setminus\{gH:g\notin H\}$.

A topological space is \textit{Polish} if it is separable and completely metrizable. For further details in descriptive set theory, we refer to~\cite{kechris}.
We say a topological group is \textit{Polish} if its topology is Polish. Consider a Polish group $G$ and a Polish space $X$. A continuous action of $G$ on $X$,
denoted by $G\curvearrowright X$, is a continuous map $a:G\times X\to X$ which satisfies that $a(1_G,x)=x$ and $a(gh,x)=a(g,a(h,x))$ for $g,h\in G$ and $x\in X$.
For brevity, we write $gx$ in place of $a(g,x)$. The \textit{orbit equivalence relation} $E_G^X$ is defined as
$$xE_G^Xy\iff\exists g\in G\,(gx=y).$$

A Polish group is \textit{non-archimedean} if it has a neighborhood basis of the identity element consisting of open subgroups.
A metric $d$ on a group $G$ is \textit{left-invariant} if $d(gh,gk)=d(h,k)$ for all $g,h,k\in G$; we also define \textit{right-invariant} metric similarly.
We say that $d$ is \textit{two-sided invariant} if it is both left-invariant and right-invariant.
The Birkhoff-Kakutani theorem asserts that every metrizable topological group admits a compatible left-invariant metric (cf.~\cite[Theorem 2.1.1]{gaobook_1}).
We say a Polish group $G$ is CLI if it admits a compatible complete left-invariant metric; and say $G$ is TSI if it admits a compatible two-sided invariant metric.  A compatible two-sided invariant metric on a Polish group is necessarily complete (cf.~\cite[Corollary 2.2.2]{gaobook_1}).
Clearly, every TSI Polish group is also CLI. All compact or abelian Polish groups are TSI, and all locally compact Polish groups are CLI (see exercises 2.1.8 and 2.2.5 of~\cite{gaobook_1}).

Given two equivalence relations $E$ and $F$ on two sets $X$ and $Y$ respectively, we say a map $f:X\rightarrow Y$ is an \textit{$(E,F)$-homomorphism} if
$$xEy\Rightarrow f(x)Ff(y)$$
for any $x,y\in X$. Moreover, $f$ is called a \textit{reduction} of $E$ to $F$ if
$$xEy\Longleftrightarrow f(x)F f(y)$$
for all $x,y\in X$. In particular, if both $X,Y$ are Polish spaces, we say $E$ is \textit{Borel reducible} to $F$, denoted by $E\leq_B F$, if there exists a Borel reduction of $E$ to $F$.
We write $E\sim_B F$ if both $E\leq_B F$ and $F\leq_B E$ hold; and write $E<_B F$ if $E\leq_B F$ and $F\nleq_B E$. We refer to~\cite{gaobook_1} for background on Borel reducibility.

We recall some benchmark equivalence relations in the research of Borel reducibility.  The equivalence relation $E_0$ on $2^\omega$ is defined as
$$xE_0y\iff\exists m\,\forall n>m\,(x(n)=y(n)).$$
If $E$ is an equivalence relation on a Polish space $X$, we define an equivalence relation $E^\omega$ on $X^\omega$ as
$$xE^\omega y\iff\forall n\,(x(n)Ey(n)).$$
The equivalence relation $\R^\omega/c_0$ on $\R^\omega$ is defined as
$$x\R^\omega/c_0 y \iff \lim_nx(n)-y(n)=0. $$

Now we recall the definition of equivalence relations induced by Polish groups as below, and list some relevant notions and results that will be repeatedly used throughout the rest of this article.

\begin{definition}[{\cite[Definition 3.1]{DZ}}]
Let $G$ be a Polish group. We define an equivalence relation $E(G)$ on $G^\omega$ as: for $x,y\in G^\omega$,
$$xE(G)y\iff\lim_nx(n)y(n)^{-1}\mbox{ converges}.$$
We say $E(G)$ is the \textit{equivalence relation induced by $G$}. Moreover, we define $c(G)=\{x\in G^\omega:\lim_nx(n)\mbox{ converges}\}$.
Then we have
$$xE(G)y\iff xy^{-1}\in c(G)\iff c(G)x=c(G)y.$$
\end{definition}

For TSI Polish groups $G$, it is more convenient to take the following $E_*(G)$ as research object than $E(G)$.

\begin{definition}[{\cite[Definition 6.1]{DZ}}]
Let $G$ be a Polish group. We define an equivalence relation $E_*(G)$ on $G^\omega$ as: for $x,y \in G^\omega$,
$$x E_*(G)y\iff\lim_nx(0)x(1)\dots x(n)y(n)^{-1}\dots y(1)^{-1}y(0)^{-1}\mbox{ converges}.$$
\end{definition}

It is trivial that $E(G)\sim_B E_*(G)$ (cf.~\cite[Proposition 2.2]{DZlcoabe}). In this article, we will use these two equivalence relations interchangeably without further explanation.

For brevity, we define
$$(x,y)|_m^n=x(m)\cdots x(n)y(n)^{-1}y(m)^{-1}$$
for $x,y\in G^\omega$ and $m\leq n$. It is clear that
$$x E_*(G)y\iff\lim_n(x,y)|_0^n\mbox{ converges.}$$

Let $G$ be a TSI Polish group and $d$ a compatible complete two-sided invariant metric on $G$. Define for the sake of brevity
$$d(x,y)|_{[m,n+1)}=d(x,y)|_{[m,n]}=d(x(m)\dots x(n),y(m)\dots y(n))$$
for $x,y\in G^\omega$ and $m\leq n$. Clearly, we have $d(x,y)|_{[m,n]}=d(1_G,(x,y)|_m^n)$.

\begin{proposition}[{\cite[Proposition 3.4]{DZ}}]\label{pro closed group}
Let $G,H$ be two Polish groups. If $G$ is topologically isomorphic to a closed subgroup of $H$, then $E(G)\leq_B E(H)$.
\end{proposition}

\begin{lemma}[{\cite[Lemma 6.2]{DZ}}]\label{Lemm basic Tsi property}
Let $G$ be a TSI Polish group, $d$ a compatible complete two-sided invariant metric on $G$.
\begin{enumerate}
\item[(1)] For $g_0,\cdots,g_n,h_0,\cdots,h_n\in G$, we have
$$d(g_0\cdots g_n,h_0\cdots h_n)=d(g_0\cdots g_nh_n^{-1}\cdots h_0^{-1},1_G)\le\sum_{k=0}^nd(g_k,h_k).$$
\item[(2)] For $x,y\in G^\omega$, we have
$$x E_*(G)y\iff\lim_m\sup_{n>m}{d(x,y)|_{[m,n]}}=0.$$
\item[(3)] For $x,y\in G^\omega$, if $xE_*(G)y$, then $\lim_nd(x(n),y(n))=0$.
\end{enumerate}
\end{lemma}

Now we recall the notion of additive reduction and a powerful lemma which converts a Borel reduction to an additive reduction.

\begin{definition}[{Farah~\cite{farah}}]
\begin{enumerate}
\item[(1)] A map $\psi:\prod_nX_n\to\prod_nX_n'$ is \textit{additive} if there exist $0=l_0<l_1<\cdots<l_j<\cdots$ and maps $T_j:X_j\to\prod_{n\in[l_j,l_{j+1})}X_n'$ such that, for $x\in\prod_nX_n$,
    $$\psi(x)=T_0(x(0))^\smallfrown T_1(x(1))^\smallfrown T_2(x(2))^\smallfrown\cdots.$$
\item[(2)] Let $E$ and $F$ be two equivalence relations on $\prod_nX_n$ and $\prod_nX_n'$ respectively, we say that $E$ is \textit{additive reducible} to $F$, denote by $E\le_AF$, if there exists an additive reduction of $E$ to $F$.
\end{enumerate}
\end{definition}

Let $E$ be an equivalence relation on $\prod_nX_n$, and let $I\subseteq\omega$ be infinite. Fix an element $w\in\prod_{n\notin I}X_n$.
For $x\in\prod_{n\in I}X_n$, define $x\oplus w\in\prod_nX_n$ as: $(x\oplus w)(n)=x(n)$ for $n\in I$ and $(x\oplus w)(n)=w(n)$ for $n\notin I$.
We define $E|_I^w$ on $\prod_{n\in I}X_n$ as: for $x,y\in\prod_{n\in I}X_n$,
$$xE|_I^wy\iff(x\oplus w)E(y\oplus w).$$

Let $(F_n)$ be a sequence of finite sets. A special equivalence relation $E_0(\prod_nF_n)$ defined as: for $x,y\in\prod_nF_n$,
$$xE_0(\prod_nF_n)y\iff\exists m\,\forall n>m\,(x(n)=y(n)).$$

\begin{lemma}[{\cite[Lemma 6.9]{DZ}}]\label{Borel_to_Additive}
Let $E$ be a Borel equivalence relation on $\prod_nF_n$ with $E_0(\prod_nF_n)\subseteq E$, where all $F_n$ are finite sets. Let $H$ be a TSI Polish group. If $E\le_BE_*(H)$, then there exist an infinite set $I\subseteq\omega$ and a $w\in\prod_{n\notin I}F_n$ such that $E|_I^w\le_AE_*(H)$.
\end{lemma}

The following notions and results are not directly presented in~\cite{DZ}, but ideas of them have already appeared in it.

\begin{definition}\label{omega and sharp}
Given two sets $X,Y$ and a map $S:X\to Y$, we define two maps $S^\omega:X^\omega\to Y^\omega$ and $S^\#:X^\omega\to Y^\omega\times X^\omega$ as: for $x\in X^\omega$ and $n\in\omega$,
$$S^\omega(x)(n)=S(x(n)),$$
$$S^\#(x)=(S^\omega(x),x).$$
\end{definition}

For any metric space $(M,d)$, recall that $E(M;0)$ is an equivalence relation on $M^\omega$ (cf.~\cite[Definition 3.2]{ding12}) defined as: for $x,y\in M^\omega$,
$$xE(M;0)y\iff\lim_nd(x(n),y(n))=0.$$
Note that, for any Polish group $G$, the equivalence relation $E(G;0)$ is independent of any choice of left-invariant compatible metric $d$ on $G$, since $d(x(n),y(n))\to 0$ iff $x(n)^{-1}y(n)\to 1_G$.

\begin{proposition}\label{equiv for sharp}
Let $G,H$ be two TSI Polish groups, and let $S:G\to H$. Then the following are equivalent:
\begin{enumerate}
\item[(1)] If $\lim_nx(n)^{-1}y(n)=1_G$, then it holds for all $x,y\in G^\omega$ that
  $$xE_*(G)y\iff S^\omega(x)E_*(H)S^\omega(y).$$
\item[(2)] $S^\#$ is a reduction of $E_*(G)$ to $E_*(H)\times E(G;0)$.
\end{enumerate}
\end{proposition}

\begin{proof}
This is an easy corollary of Lemma~\ref{Lemm basic Tsi property}(3).
\end{proof}

\begin{lemma}\label{converse}
Let $G,H$ be two TSI Polish groups such that $E_*(G)\le_B E_*(H)\times E(G;0)$. If the interval $[0,1]$ can be embedded into $H$, or $E(G;0)\le_B E(H;0)$ holds, then we have $E_*(G)\le_B E_*(H)$.
\end{lemma}

\begin{proof}
By~\cite[Theorem 3.4(ii)]{ding12}, we have $E(G;0)\le_B E([0,1];0)$.

Let $f:[0,1]\to H$ be an embedding. By the uniformly continuity of $f$ and $f^{-1}:f([0,1])\to[0,1]$, we see that $E([0,1];0)\sim_B E(f([0,1]);0)\le_B E(H;0)$. So $E_*(G)\le_B E_*(H)\times E(H;0)$.

Finally, by~\cite[Lemma 6.3]{DZ}, we obtain $E_*(H)\times E(H;0)\le_B E_*(H)$.
\end{proof}

\section{Non-TSI vs TSI Polish groups}

In this section, we show that we can use Borel reducibility among equivalence relations induced by Polish groups to characterize the class of TSI Polish groups.

Let $G$ be a CLI Polish group, $H$ a TSI Polish group. Suppose that $d_G'$ is a compatible complete left-invariant metric on $G$ and $d_H$ is a compatible complete two-sided invariant metric on $H$. Put $d_G(g,g')=d_G'(g^{-1},(g')^{-1})$ for $g,g'\in G$. It is trivial to check that $d_G$ is a compatible complete right-invariant metric on $G$.

Assume that $E_*(G)\le_BE_*(H)$. Let $(F_n)$ be a sequence of finite subsets of $G$ such that
\begin{enumerate}
\item[(i)] $1_G\in F_n=F_n^{-1}$;
\item[(ii)] $F_{n-1}^{n}\subseteq F_n$; and
\item[(iii)] $\bigcup_nF_n$ is dense in $G$.
\end{enumerate}
Denote by $E$ the restriction of $E_*(G)$ on $\prod_nF_n$. By Lemma~\ref{Borel_to_Additive}, there exist an infinite set $I\subseteq\omega$ and a $w\in\prod_{n\notin I}F_n$ such that $E|_I^w\le_AE_*(H)$. So there are natural numbers $0=n_0<n_1<n_2<\cdots$ with $I=\{n_j:j\in\omega\}$, $0=l_0<l_1<l_2<\cdots$, maps $T_{n_j}:F_{n_j}\to H^{l_{j+1}-l_j}$, and $\psi:\prod_{n\in I}F_n\to H^\omega$ with
$$\psi(x)=T_{n_0}(x(n_0))^\smallfrown T_{n_1}(x(n_1))^\smallfrown T_{n_2}(x(n_2))^\smallfrown\cdots,$$
such that $\psi$ is an additive reduction of $E|_I^w$ to $E_*(H)$.

Put $u_{n_0}=1_G$, and for $j>0$, let
$$u_{n_j}=w(n_{j-1}+1)\cdots w(n_j-1).$$
By (i) and (ii), $u_{n_j}^{-1}\in F_{n_j}$ for each $j\in\omega$.

Define $x_0\in\prod_{n\in I}F_n$ as
$$x_0(n_j)=u_{n_j}^{-1}\quad(\forall j\in\omega).$$
We may assume that $\psi(x_ 0)=1_{H^\omega}$, otherwise we can replace $\psi$ with the following $\psi'$: for $x\in\prod_{n\in I}F_n$ and $k\in\omega$,
$$\psi'(x)(k)=\psi(x_0)(0)\cdots\psi(x_0)(k-1)\psi(x)(k)\psi(x_0)(k)^{-1}\cdots\psi(x_0)(0)^{-1}.$$
Clearly, $(\psi(x),\psi(y))|_0^k=(\psi'(x),\psi'(y))|_0^k$, so we have
$$\psi(x)E_*(H)\psi(y)\iff\psi'(x)E_*(H)\psi'(y)$$
holds for $x,y\in\prod_{n\in I}F_n$. Hence $\psi'$ is an additive reduction of $E|_I^w$ to $E_*(H)$ with $\psi'(x_0)=1_{H^\omega}$, as desired.

For $s=(h_0,\ldots,h_{l-1})$ and $t=(h_0',\ldots,h_{l-1}')$ in $H^l$, we define
$$d_H^\infty(s,t)=\max_{0\le i\le m<l}d_H(h_i\cdots h_m,h_i'\cdots h_m').$$

For any integer $j>0$ and $g\in F_{n_j-1}$, it is worth noting that
$$u_{n_j}^{-1}g=w(n_j-1)^{-1}\cdots w(n_{j-1}+1)^{-1}g\in F_{n_j-1}^{n_j}\subseteq F_{n_j}.$$

\begin{lemma}\label{continuous property}
For any $q\in\omega$, there exists a $\delta_q>0$ such that
$$\forall^\infty n\in I\,\forall g,g'\in F_{n-1}\,(d_G(g,g')<\delta_q\Rightarrow d_H^\infty(T_n(u_n^{-1}g),T_n(u_n^{-1}g'))<2^{-q}).$$
\end{lemma}

\begin{proof}
If not, then there exist an $\varepsilon_0>0$, a strictly increasing sequence $(j_k)$, and $g_k,g_k'\in F_{n_{j_k}-1}$ for each $k$, such that
$$d_G(g_k,g_k')<2^{-(k+1)},\quad d_H^\infty(T_{n_{j_k}}(u_{n_{j_k}}^{-1}g_k),T_{n_{j_k}}(u_{n_{j_k}}^{-1}g_k'))\ge\varepsilon_0.$$

Since $d_G$ is right-invariant, we have $\lim_kg_k'g_k^{-1}=1_G$. We shall inductively find a strictly increasing sequence of natural numbers $k_0<k_1<\cdots$ so that
$$d_G(g_{k_0}'\cdots g_{k_{p-1}}'g_{k_{p-1}}^{-1}\cdots g_{k_0}^{-1},g_{k_0}'\cdots g_{k_p}'g_{k_p}^{-1}\cdots g_{k_0}^{-1})<2^{-p}$$
for each integer $p>0$. Put $k_0=0$ first. Assume that $k_0,\ldots,k_{p-1}$ have been found. By the continuity of group operations, there is some $\delta>0$ such that, for $g,g'\in G$ with $d_G(1_G,g'g^{-1})<\delta$, we have
$$d_G(g_{k_0}'\cdots g_{k_{p-1}}'1_Gg_{k_{p-1}}^{-1}\cdots g_{k_0}^{-1},g_{k_0}'\cdots g_{k_{p-1}}'g'g^{-1}g_{k_{p-1}}^{-1}\cdots g_{k_0}^{-1})<2^{-p}.$$
Then we can find a $k_p>k_{p-1}$ large enough so that $d_G(1_G,g_{k_p}'g_{k_p}^{-1})<\delta$. This completes the induction. It follows that, for $m>p$, we have
$$d_G(g_{k_0}'\cdots g_{k_p}'g_{k_p}^{-1}\cdots g_{k_0}^{-1},g_{k_0}'\cdots g_{k_m}'g_{k_m}^{-1}\cdots g_{k_0}^{-1})<\sum_{i=p+1}^m 2^{-i}<2^{-p}.$$
Since $d_G$ is complete, we get that
$$\lim_pg_{k_0}'\cdots g_{k_p}'g_{k_p}^{-1}\cdots g_{k_0}^{-1}\mbox{ converges.}$$

For each $n\in I$, put
$$x(n)=\left\{\begin{array}{ll}u_n^{-1}g_{k_p}', & n=n_{j_{k_p}},\cr u_n^{-1}, & \mbox{otherwise,}\end{array}\right.\quad
y(n)=\left\{\begin{array}{ll}u_n^{-1}g_{k_p}, & n=n_{j_{k_p}},\cr u_n^{-1}, & \mbox{otherwise.}\end{array}\right.$$
For any $n\ge n_{j_{k_0}}$, let $p_n$ be the largest $p$ such that $n_{j_{k_p}}\le n$. Note that
$$(x\oplus w,y\oplus w)|_0^n=g_{k_0}'\cdots g_{k_{p_n}}'g_{k_{p_n}}^{-1}\cdots g_{k_0}^{-1},$$
so $(x\oplus w)E_*(G)(y\oplus w)$. But for each $p\in\omega$, if $k=k_p$, then
$$\max_{l_{j_k}\le i\le m<l_{j_k+1}}d_H(\psi(x),\psi(y))|_{[i,m]}=
d_H^\infty(T_{n_{j_k}}(u_{n_{j_k}}^{-1}g_k'),T_{n_{j_k}}(u_{n_{j_k}}^{-1}g_k))\ge\varepsilon_0.
$$
So $\psi(x)E_*(H)\psi(y)$ fails, contradicting that $\psi$ is a reduction.
\end{proof}

\begin{definition}\label{S_n_j}
For any $g\in\bigcup_nF_n$, if $g\in F_n$ for some $n<n_j$, then $u_{n_j}^{-1}g\in F_{n_j}$, so $T_{n_j}(u_{n_j}^{-1}g)\in H^{l_{j+1}-l_j}$. This allows us to define
$$S_{n_j}(g)=T_{n_j}(u_{n_j}^{-1}g)(0)\cdots T_{n_j}(u_{n_j}^{-1}g)(l_{j+1}-l_j-1)\in H.$$
\end{definition}

Recall that $\psi(x_ 0)=1_{H^\omega}$, this implies that
$$S_{n_j}(1_G)=T_{n_j}(u_{n_j}^{-1})(0)\cdots T_{n_j}(u_{n_j}^{-1})(l_{j+1}-l_j-1)=1_H$$
for all $j\in\omega$. It is clear that, for $g,g'\in F_n$ with $n<n_j$,
$$d_H(S_{n_j}(g),S_{n_j}(g'))\le d_H^\infty(T_{n_j}(u_{n_j}^{-1}g),T_{n_j}(u_{n_j}^{-1}g')).$$

\begin{lemma}\label{preserving converges on F_n}
Let $x,y\in\prod_{n\in I}F_n$ such that $\lim_j d_G(x(n_j),y(n_j))=0$ and $x(n_j),y(n_j)\in F_{n_j-1}$ for all $j>0$. Then
$$\psi(x') E_*(H)\psi(y')\iff\lim_j(\psi(x'),\psi(y'))|_0^{l_j-1}\rm{ converges},$$
where $x'(n)=u_{n}^{-1}x(n),y'(n)=u_{n}^{-1}y(n)$ for all $n\in I$.
\end{lemma}

\begin{proof}
If $\psi(x') E_*(H)\psi(y')$,
then $\lim_p(\psi(x'),\psi(y'))|_0^{p}\mbox{ converges}.$
In particular,  $\lim_j(\psi( x'),\psi(y'))|_0^{l_j-1}\mbox{converges}.$

Conversely, suppose $\lim_j(\psi(x'),\psi(y'))|_0^{l_j-1}$ converges. For any $q\in\omega$, by Lemma~\ref{continuous property}, there is a $\delta_q\in\omega$ such that
$$\forall^\infty n\in I\,\forall g,g'\in F_{n-1}\,(d_G(g,g')<\delta_q\Rightarrow d_H^\infty(T_n(u_n^{-1}g),T_n(u_n^{-1}g'))<2^{-q}).$$
Since $\lim_j d_G(x(n_j),y(n_j))=0$ and $x(n_j),y(n_j)\in F_{n_j-1}$ for any integer $j>0$, there exists a $j_0\in\omega$ such that
$$\forall j>j_0\, (d_H^\infty(T_{n_j}(u_{n_j}^{-1}x(n_j)),T_{n_j}(u_{n_j}^{-1}y(n_j)))<2^{-q}).$$
Then for any $j,k\in\omega$ with $j>j_0$ and $l_j\leq k<l_{j+1}$, since $d_H$ is two-sided invariant, we have
$$\begin{array}{ll} d_H((\psi(x'),\psi(y'))|_0^{l_j-1},(\psi(x'),\psi(y'))|_0^{k})&=d_H(\psi(x'),\psi(y'))|_{[l_j,k]}\cr
&\leq d_H^\infty(T_{n_j}(x'(n_j)),T_{n_j}(y'(n_j)))\cr
&<2^{-q}.\end{array}$$
So by the convergency of $\lim_j(\psi(x'),\psi(y'))|_0^{l_j-1}$, we obtain that
$$\lim_k(\psi(x'),\psi(y'))|_0^{k}\mbox{ converges},$$
and hence $\psi(x') E_*(H)\psi(y')$.
\end{proof}

The following proposition may be well known. We provide a proof here for the convenience of readers.

\begin{proposition}[folklore]\label{prop for TSI}
Let $G$ be a Polish group. Then $G$ is TSI iff for all sequences $(g_p)$ and $(g_p')$ in $G$, we have
$$\lim_pg_pg_p'=1_G\iff\lim_pg_p'g_p=1_G.$$
\end{proposition}

\begin{proof}
Suppose $G$ is TSI, and let $d_G$ be a compatible complete two-sided invariant metric on $G$. Then we have $d_G(1_G,g_pg_p')=d_G(g_p^{-1},g_p')=d_G(1_G,g_p'g_p)$, so $\lim_pg_pg_p'=1_G\iff\lim_pg_p'g_p=1_G$.

On the other hand, fix an open neighborhood basis $(V_n)$ of $1_G$. We claim that
$$\forall n\,\exists m_n\,\forall g\in G\,(gV_{m_n}g^{-1}\subseteq V_n).$$
If not, there exist an $n_0$ and two sequences $(g_p),(h_p)$ in $G$ such that $\lim_ph_p=1_G$, but $g_ph_pg_p^{-1}\notin V_{n_0}$ for each $p\in\omega$. Put $g_p'=h_pg_p^{-1}$. It is clear that $\lim_pg_p'g_p=1_G$, but $\lim_pg_pg_p'\ne 1_G$, which is a contradiction.

Now we put $U_n=\bigcup_{g\in G}gV_{m_n}g^{-1}$. Note that $U_n\subseteq V_n$ for each $n\in\omega$. Therefore, $(U_n)$ is also an open neighborhood basis of $1_G$ with $gU_ng^{-1}=U_n$ for all $g\in G$ and $n\in\omega$. By~\cite[Exercise 2.1.4]{gaobook_1}, $G$ admits a compatible two-sided invariant metric. Hence $G$ is TSI (cf.~\cite[Corollary 1.2.2]{BK}).
\end{proof}

\begin{theorem}\label{TSI}
Let $G,H$ be two Polish groups. If $H$ is TSI but $G$ is not, then $E(G)\not\le_BE(H)$.
\end{theorem}

\begin{proof}
Assume toward a contradiction that $E(G)\leq_B E(H)$. By Theorem 4.3 of~\cite{DZ}, it suffices to consider the case that $G$ is CLI. Let $d_G$ be a compatible complete right-invariant metric on $G$ and $d_H$ a compatible complete two-sided invariant metric on $H$.

We use the notation defined earlier in this section.

Since $G$ is not TSI, by Proposition~\ref{prop for TSI}, there are two sequences $(g_p),(g_p')$ of elements of $G$ such that $\lim_pg'_pg_p=1_G$, but $\lim_pg_pg_p'\ne 1_G$. Since $d_G$ is right-invariant, $\lim_p d_G(g'_p,g_p^{-1})=0$. By transfering to a subsequence, we may assume that, there is a $\delta>0$ such that $\inf_p d_G(1_G,g_pg'_p)>\delta$.

Since $\bigcup_n F_n$ is dense, by perturbation, we may assume that $\{g_p,g_p':p\in\omega\}\subseteq\bigcup_n F_n$. By Lemma~\ref{continuous property} and Definition~\ref{S_n_j}, we can find two sequence $(p(i))$ and $(q(i))$ of natural numbers so that
\begin{enumerate}
   \item[(i)] $g_{p(i)},g_{p(i)}'\in F_{n_{q(i)}-1}$;
   \item[(ii)] $d_H(S_{n_{q(i)+1}}(g'_{p(i)}),S_{n_{q(i)+1}}(g_{p(i)}^{-1}))<2^{-i}$; and
   \item[(iii)] for each $i\in\omega$, $0<p(i)<q(i)<q(i)+1<p(i+1)$.
\end{enumerate}
For each $n\in I$, define
$$x(n)=\left\{\begin{array}{ll}
   u_n^{-1} g_{p(i)}, & n=n_{q(i)}, \cr
   u_n^{-1}g'_{p(i)}, & n=n_{q(i)+1},\cr
   u_n^{-1},& \mbox{otherwise},
  \end{array}\right. \quad
  y(n)=\left\{\begin{array}{ll}
   u_n^{-1}g_{p(i)}, & n=n_{q(i)}, \cr
   u_n^{-1}g_{p(i)}^{-1}, & n=n_{q(i)+1},\cr
   u_n^{-1},& \mbox{otherwise}.
  \end{array}\right.
  $$
For $j\in\omega$, by letting $h'_j=S_{n_j}(u_{n_j}x(n_j))$ and $h_j=S_{n_j}(u_{n_j}y(n_j))$, we have
$$(\psi(x),\psi(y))|_0^{l_{j+1}-1}=h'_0\cdots h'_jh_j^{-1}\cdots h_0^{-1}.$$
For $i\in\omega$ and $m>j>q(i)+1$, it follows from (ii) and Lemma~\ref{Lemm basic Tsi property}(1) that
$$d_H(h'_j\cdots h'_m,h_j\cdots h_m)<\sum_{k>i} 2^{-k}<2^{-i}.$$
Now by Lemma~\ref{Lemm basic Tsi property}(2),
$$\lim_j(\psi(x),\psi(y))|_0^{l_{j+1}-1}=\lim_jh'_0\cdots h'_jh_j^{-1}\cdots h_0^{-1}\mbox{ converges}.$$
Note that $\lim_i d_G(g'_{p(i)},g_{p(i)}^{-1})=0$. So by Lemma~\ref{preserving converges on F_n}, we have
$$\psi(x) E_*(H)\psi(y),$$
and hence $(x\oplus w) E_*(G)(y\oplus w)$. In particular,
$$\lim_i(x\oplus w,y\oplus w)|_0^{n_{q(i)+1}}=\lim_i g_{p(0)}g'_{p(0)}\cdots g_{p(i)}g'_{p(i)}\mbox{ converges}.$$
Therefore, it follows that
$$\lim_i g_{p(i)}g'_{p(i)}=\lim_i(g_{p(0)}g'_{p(0)}\cdots g_{p(i-1)}g'_{p(i-1)})^{-1}g_{p(0)}g'_{p(0)}\cdots g_{p(i)}g'_{p(i)}=1_G.$$
We obtain a contradiction.
\end{proof}

\section{A gap between $E_0^\omega$ and $\R^\omega/c_0$}

In this section, we prove that there is NO Polish group $G$ satisfying that
$$E_0^\omega<_B E(G)\leq_B \R^\omega/ c_0.$$
To do so, we will use the following lemma, which is an easy corollary of~\cite[Lemma 4.2]{ding17}.

\begin{lemma}\label{Lemm Add of c0}
Let $E$ be a Borel equivalence relation on $\prod_nF_n$ with $E_0(\prod_nF_n)\subseteq E$, where all $(F_n)$ are finite sets. If $E\le_B \R^\omega/c_0$, then there exist an infinite set $I\subseteq\omega$ and a $w\in\prod_{n\notin I}F_n$ such that $E|_I^w\le_A\R^\omega/c_0$.
\end{lemma}

\begin{proof}
Let $e_n=(0,\ldots,0,\overset{n}1,0,\ldots)$ for each integer $n>0$. Then $(e_n)$ is the canonical Schauder basis of $c_0$. Note that $E(c_0,(e_n))=\R^\omega/c_0$.
By applying~\cite[Lemma 4.2]{ding17} to $E(c_0,(e_n))$, we conclude the proof.
\end{proof}

\begin{theorem}\label{theo-co}
Let $G$ be a Polish group. If $E(G)\leq_B \R^\omega/c_0$, then $G$ is TSI non-archimedean.
\end{theorem}

\begin{proof}
Suppose $E(G)\leq_B \R^\omega/c_0$. Note that $\R^\omega/c_0<_B\R^\omega/{\rm cs}=E_*(\R)$ (see~\cite[Theorem 5.9(i)]{ding17}). By Theorem~\ref{TSI}, $G$ is also TSI. Let $d_G$ be a compatible complete two-sided invariant metric on $G$. Put $V_k=\{g\in G:d_G(1_G,g)<2^{-k}\}$.

Assume for contradiction that $G$ is not non-archimedean. Then there is a $K\in\omega$ such that $V_K$ contains no open subgroups of $G$. It is clear that $V_k=V_k^{-1}$ and $\bigcup_mV_k^m$ is an open subgroup which is not contained in $V_K$. So there exists $m_k>0$ with $V_k^{m_k}\not\subseteq V_K$ for each $k$. We can find $g_{k,0},\ldots,g_{k,m_k-1}\in V_k$ such that $g_{k,0}\cdots g_{k,m_k-1}\notin V_K$. Put $h_j=g_{k,i}$ for $j=\sum_{l<k}m_l+i$ with $i<m_k$. Then we have $\lim_j h_j=1_G$. By Cauchy Criterion,
$$\lim_jh_0\cdots h_j\mbox{ diverges}.$$

Assume that $E_*(G)\le_B\R^\omega/c_0$. Let $(F_n)$ be a sequence of finite subsets of $G$ such that
\begin{enumerate}
\item[(i)] $1_G,h_0\in F_n=F_n^{-1}$;
\item[(\romannumeral 2)] $F_{n-1}^{n}\subseteq F_n$; and
\item[(\romannumeral 3)] for $n>0,h_n\in F_{n-1}$.
\end{enumerate}
Denote by $E$ the restriction of $E_*(G)$ on $\prod_nF_n$. By Lemma~\ref{Lemm Add of c0}, there exist an infinite set $I\subseteq\omega$ and a $w\in\prod_{n\notin I}F_n$ such that $E|_I^w\le_A\R^\omega/c_0$.
So there are natural numbers $0=n_0<n_1<n_2<\cdots$ with $I=\{n_j:j\in\omega\}$, $0=l_0<l_1<l_2<\cdots$, maps $T_{n_j}:F_{n_j}\to\R^{l_{j+1}-l_j}$, and $\psi:\prod_{n\in I}F_n\to\R^\omega$ with
$$\psi(x)=T_{n_0}(x(n_0))^\smallfrown T_{n_1}(x(n_1))^\smallfrown T_{n_2}(x(n_2))^\smallfrown\cdots,$$
such that $\psi$ is an additive reduction of $E|_I^w$ to $\R^\omega/c_0$.

Put $u_{n_0}=1_G$, and for $j>0$, let
$$u_{n_j}=w(n_{j-1}+1)\cdots w(n_j-1).$$

Assume again for contradiction that there are $\varepsilon_0>0$ and natural numbers $0<j(0)<j(1)<\cdots<j(p)<\cdots$ such that
$$\max_{0\leq k <l_{j(p)+1}-l_{j(p)}}|T_{n_{j(p)}}(u^{-1}_{n_{j(p)}}h_{j(p)})(k)-T_{n_{j(p)}}(u^{-1}_{n_{j(p)}})(k)|\geq\varepsilon_0.$$
By $\lim_p h_{j(p)}=1_G$, we can find natural numbers $p_0<p_1<\cdots<p_i<\cdots$ so that $d_G(1_G,h_{j(p_i)})<2^{-i}$ for each $i\in\omega$. Thus $(h_{j(p_0)}\cdots h_{j(p_i)})$ is $d_G$-Cauchy, so $\lim_ih_{j(p_0)}\cdots h_{j(p_i)}$ converges. For each $n\in I$, put $x_0(n)=u^{-1}_{n}$ and
$$y_0'(n)=\left\{\begin{array}{ll} u^{-1}_{n}h_{j(p_i)}, & n=n_{j(p_i)}, \cr u^{-1}_n, & \mbox{otherwise.}\end{array}\right.$$
For any $n> n_{j(p_0)}$, let $i_n$ be the largest $i$ with $n_{j(p_{i})}\leq n$. Then
$$(y_0'\oplus w,x_0\oplus w)|_0^n=h_{j(p_0)}\cdots h_{j(p_{i_n})},$$
and thus $(y_0'\oplus w) E_*(G)(x_0\oplus w)$. Note that for all $i\in\omega$,
$$\max_{l_{j(p_i)}\leq k<l_{j(p_i)+1}}|\psi(y_0')(k)-\psi(x_0)(k)|\geq\varepsilon_0.$$
So $\psi(y_0')\R^\omega/c_0 \psi(x_0)$ fails, which is a contradiction.

Therefore, for all $\varepsilon>0$, we have
$$\forall^\infty j\,\forall k\in[0,l_{j+1}-l_j)\,(|T_{n_j}(u^{-1}_{n_j}h_j)(k)-T_{n_j}(u^{-1}_{n_j})(k)|<\varepsilon).$$
For each $n\in I$, put $y_0(n)=u^{-1}_nh_j$, where $n=n_j$. Then we have
$$\psi(y_0)\R^\omega/c_0\psi(x_0),$$
and hence $(y_0\oplus w) E_*(G)(x_0\oplus w)$. In particular,
$$\lim_j(y_0\oplus w,x_0\oplus w)|_0^{n_j}=\lim_jh_0\cdots h_j\mbox{ converges}.$$
This leads to a contradiction.
\end{proof}

Building upon previous results, we establish the following theorem.

\begin{theorem}
Let $G$ be a Polish group. Then the following are equivalent:
\begin{enumerate}
  \item $G$ is TSI non-archimedean;
  \item $E(G)\leq_B E_0^\omega$; and
  \item $E(G)\leq_B \R^\omega/c_0$.
\end{enumerate}
In particular, $E(G)\sim_B E_0^\omega$ iff $G$ is TSI uncountable non-archimedean.
\end{theorem}

\begin{proof}
$(1)\Rightarrow(2)$ follows from~\cite[Theorem 3.5]{DZ}. $(2)\Rightarrow(3)$ follows from~\cite[Lemma 8.5.3]{gaobook_1}. $(3)\Rightarrow(1)$ follows from Theorem~\ref{theo-co}.

Again by~\cite[Theorem 3.5]{DZ}, we see that $E(G)\sim_B E_0^\omega$ iff $G$ is TSI uncountable non-archimedean.
\end{proof}

Note that $\R^\omega/c_0<_B\R^\omega/{\rm cs}=E_*(\R)$ (see~\cite[Theorem 5.9(i)]{ding17}) and $E(\R)=\R^\omega/c$.
The results of this section make the following question very interesting.

\begin{question}
For any Polish group $G$, does it hold that
$$E(G)<_BE(\R)\iff E(G)\le_BE_0^\omega?$$
\end{question}

\section{The Pre-rigid Theorem on TSI Polish groups}

In this section, we prove a highly technical theorem, namely the Pre-rigid Theorem, which will serve as the foundation for the subsequent sections.

We say that a topological group $G$ \textit{has no small subgroups}, or \textit{is NSS}, if there exists an open subset $V\ni 1_G$ in $G$ such that no non-trivial subgroup of $G$ contained in $V$.

To generalize Theorem 6.13 of~\cite{DZ}, we introduce the following definition.

\begin{definition}
A Polish group $G$ is called \textit{strongly NSS} if there exists an open set $V\ni 1_G$ in $G$ such that
$$\forall(g_n)\in G^\omega\,(\lim_n g_n\ne 1_G\Rightarrow\exists n_0<\cdots<n_k\,(g_{n_0}\cdots g_{n_k}\notin V)),$$
where the set $V$ is called an \emph{unenclosed set} of $G$.
\end{definition}

\begin{proposition}\label{NSS and strongly NSS}
Let $G$ be a Polish group. The following hold:
\begin{enumerate}
  \item [(1)] if $G$ is strongly NSS, then $G$ is NSS; and
  \item [(2)] if $G$ is locally compact, then $G$ is strongly NSS iff $G$ is NSS.
\end{enumerate}
\end{proposition}

\begin{proof}
(1) Let $V$ be an unenclosed set of $G$. We have
$$\forall g\in G\,(g\neq 1_G\Rightarrow \exists m\,(g^m\notin V)).$$
Thus $V$ contains no non-trivial subgroups of $G$, so $G$ is NSS.

(2) By (1), we only need to prove another direction. Suppose that $G$ is NSS and locally compact. Let $V\ni 1_G$ be an open subset of $G$ such that $\overline{V}$ is compact and contains no non-trivial subgroups of $G$. Then
$$\forall g\in \overline{V}\,(g\neq 1_G\Rightarrow \exists m\,(g^m\notin \overline{V})).$$

We claim that $V$ is an unenclosed set of $G$.
Fix a $(g_n)\in G^\omega$ with $\lim_n g_n\ne 1_G$. By the definition of unenclosed set, it suffices to consider the case that $(g_n)\in V^\omega$.
Since $\overline{V}$ is compact, there exist a subsequence $(g_{n_i})$ of $(g_n)$ and an element $1_G\neq h\in \overline{V}$ such that $\lim_i g_{n_i}=h$.
Then we can find an $m\in\omega$ with $h^m\notin \overline{V}$. By the continuity of group operations, there are $i_0<i_1<\cdots<i_{m-1}$ such that $g_{n_{i_0}}g_{n_{i_1}}\cdots g_{n_{i_{m-1}}}\notin \overline{V} $. Therefore $G$ is strongly NSS.
\end{proof}

The addition group $\R^\omega$ with the product topology is not strongly NSS because it is not NSS.  Let $(e_n)$ be the canonical basis in $c_0$. For any $r>0$, the sequence $(\frac{r}{2}e_n)$ witnesses that the open set $\{x\in c_0:\Vert x\Vert<r\}$ is not an unenclosed set of $c_0$. On the other hand, for any $x\in c_0$ with $\Vert x\Vert\neq 0$, we have $\lim_n\Vert nx\Vert\rightarrow +\infty$. Thus the Banach space $c_0$ is NSS, but not strongly NSS, under the addition operation.

Now, we are ready to prove the Pre-rigid Theorem. Let us remind that the definitions of the maps $(\pi^S_m)^\omega$ and $S^\#$ appeared in the following theorem can be found in Definition~\ref{omega and sharp}.

\begin{theorem}[Pre-rigid Theorem]\label{Pre-rigid Theorem}
Let $G$ be a TSI Polish group, and let $H$ be a closed subgroup of the product of a sequence of TSI strongly NSS Polish groups $(H_m)$.
If $E(G)\le_BE(H)$, then for each $m\in\omega$, there exist an open subgroup $W_m$ of $G$ and a continuous map $\pi^S_m:W_m\to H_m$ with $\pi^S_m(1_{G})=1_{H_m}$ satisfying the following:
\begin{enumerate}
\item [(i)] The map $(\pi^S_m)^\omega:W_m^\omega\to H_m^\omega$ is an $(E_*(W_m),E_*(H_m))$-homomorphism.
\item [(ii)] Define $S:W=\bigcap_mW_m\to H$ as: for $g\in W$,
  $$S(g)=(\pi^S_0(g),\pi^S_1(g),\ldots,\pi^S_m(g),\ldots),$$
  then the map $S^\#:W^\omega\to H^\omega\times W^\omega$ is a reduction of $E_*(W)$ to $E_*(H)\times E(W;0)$.
\end{enumerate}

In particular, the converse holds if $G=W$ and the interval $[0,1]$ can be embedded into $H$.
\end{theorem}

\begin{proof}
Let $d_G,d_{H_m}\leq 1$ be compatible complete two-sided invariant metrics on $G$ and $H_m$ respectively.
A compatible complete two-sided invariant metric $d_H$ on $H$ is defined as
$$d_H(h,h')=\sum_{m=0}^\infty 2^{-m}d_{H_m}(h(m),h'(m))$$
for $h,h'\in H$. Assume that $E_*(G)\le_BE_*(H)$.

We use the notation defined in the arguments before Lemma~\ref{continuous property} and before Lemma~\ref{preserving converges on F_n}.

For $m\in\omega$, let $\pi_m:\prod_n H_n\to H_m$ be the canonical projection map, so $\pi_m(h)=h(m)$ for each $h$.
It is clear that every $\pi_m$ is a continuous homomorphism.

It is worth noting that
$$\forall h,h'\in H\, (d_{H_m}(\pi_m(h),\pi_m(h'))\leq 2^md_H(h,h')).$$

For $m\in\omega$, let $X_m$ be the set of all $g\in \bigcup_n F_n$ such that:
$$\forall\delta>0\, \exists j_\delta\,(g\in F_{n_{j_\delta}-1}\wedge \forall j'> j\geq j_\delta\,(d_{H_m}(1_{H_m},\pi_m(S_{n_j}(g)S_{n_{j'}}(g^{-1})))<\delta)).$$
Note that $S_{n_j}(1_G)=1_H$ for all $j\in\omega$, so $1_G\in X_m$.

\begin{claim}\label{claim conver g g-}
For any $m\in\omega$, there is an $r_m>0$ such that
$$\forall g,g'\in\bigcup_n F_n\,((d_G(g,g')<r_m\wedge g\in X_m)\Rightarrow g'\in X_m).$$
\end{claim}

\noindent\emph{Proof of Claim 1.}
If not, there are an $m_0\in\omega$ and two sequences $(g_i),(g'_i)\in\bigcup_n F_n$ such that
\begin{enumerate}
  \item [(i)] $d_G(g_i,g'_i)<2^{-i}$; and
  \item [(ii)] $g_i\in X_{m_0}$ and $g'_i\notin X_{m_0}$.
\end{enumerate}

Note that $H_{m_0}$ is strongly NSS. Let open set $V\ni 1_{H_{m_0}}$ be an unenclosed set of $H_{m_0}$. Put $p(-1,2M_{-1}-1)=0$. For any $i\in\omega$, we will find an $M_i\in\omega$ and natural numbers $0<p(i,0)<p(i,1)<\cdots<p(i,2M_i-1)$ so that
\begin{enumerate}
  \item [(1)] $p(i-1,2M_{i-1}-1)<p(i,0)$;
  \item [(2)]  $g_i,g_i'\in F_{n_{p(i,0)-1}}$;
  \item [(3)] for $i\in\omega$,
  $$\pi_{m_0}(S_{n_{p(i,0)}}(g'_i)S_{n_{p(i,1)}}((g'_i)^{-1})\cdots S_{n_{p(i,2M_i-2)}}(g'_i)S_{n_{p(i,2M_i-1)}}((g'_i)^{-1}))\notin V;$$
   \item [(4)] for $i\in\omega$, $$\sum_{l=0}^{M_i-1}d_{H_{m_0}}(1_{H_{m_0}},\pi_{m_0}(S_{n_{p(i,2l)}}(g_i)S_{n_{p(i,2l+1)}}({g_i^{-1}})))<2^{-i}.$$
\end{enumerate}

Let us begin with $i=0$. Since $g'_0\notin X_{m_0}$, there are a $\delta_0>0$ and natural numbers $p(0)<p(1)<\cdots<p(j)<\cdots$ such that $g'_0\in F_{n_{p(0)}-1}$ and for each $j\in\omega$,
$$d_{H_{m_0}}(1_{H_{m_0}},\pi_{m_0}(S_{n_{p(2j)}}(g'_0)S_{n_{p(2j+1)}}((g_0')^{-1})))\geq\delta_0.$$
By the definition of $X_{m_0}$ and $g_0\in X_{m_0}$, there exists a strictly increasing sequence $(j_k)$ of natural numbers such that $g_0\in F_{n_{p(j_0)}-1}$ and for each $k\in\omega$,
$$d_{H_{m_0}}(1_{H_{m_0}},\pi_{m_0}(S_{n_{p(2j_k)}}(g_0)S_{n_{p(2j_k+1)}}({g_0^{-1}})))<2^{-(k+1)}.$$
Let $h\in H_{m_0}^\omega $ be such that $$h(k)=\pi_{m_0}(S_{n_{p(2j_{k})}}(g'_0)S_{n_{p(2j_{k}+1)}}((g'_0)^{-1}))$$
for each $k\in\omega$. Then $\lim_k h(k)\ne 1_{H_{m_0}}$. Since $V$ is an unenclosed set of $H_{m_0}$, there are finitely many natural numbers $k_0<k_1<\cdots<k_{q}$ such that $h({k_0})\cdots h({k_{q}})\notin V$, i.e.,
$$\pi_{m_0}(S_{n_{p(2j_{k_0})}}(g'_0)S_{n_{p(2j_{k_0}+1)}}((g'_0)^{-1})\cdots S_{n_{p(2j_{k_{q}})}}(g'_0)S_{n_{p(2j_{k_{q}}+1)}}((g'_0)^{-1}))\notin V.$$
Now we put $M_0=q+1$, $p(0,2l)=p(2j_{k_l})$, and $p(0,2l+1)=p(2j_{k_l}+1)$ for all $l<q+1$. Then we have $g_0,g_0'\in F_{n_{p(j_0)}-1}\subseteq F_{n_{p(0,0)}-1}$, and
$$\sum_{l=0}^{q}d_{H_{m_0}}(1_{H_{m_0}},\pi_{m_0}(S_{n_{p(0,2l)}}(g_0)S_{n_{p(0,2l+1)}}({g_0^{-1}})))<\sum_{l=0}^{q}2^{-(k_l+1)}<1.$$
We can see that clauses (1)--(4) hold for $i=0$.

Now assume that $M_0,\ldots,M_{i}$ and $p(0,0)<\cdots<p(0,2M_0-1)<\cdots<p(i,0)<\cdots<p(i,2M_{i}-1)$ have been defined. Using similar arguments as those presented in the preceding paragraph, pick a $p'>p(i,2M_{i}-1)$ with $g_{i+1},g_{i+1}'\in F_{n_{p'}-1}$, then we can find an $M_{i+1}\in\omega$ and natural numbers $p'\le p(i+1,0)<\cdots<p(i+1,2M_{i+1}-1)$ so that clauses (1)--(4) hold.

For each $n\in I$, define
$$x(n)=\left\{\begin{array}{ll}
   u_n^{-1} g'_{i}, &  n=n_{p(i,2k)},0\leq k<M_i , \cr
   u_n^{-1}{(g_i')}^{-1}, & n=n_{p(i,2k+1)},0\leq k<M_i ,\cr
   u_n^{-1},& \mbox{otherwise},
  \end{array}\right.$$

$$y(n)=\left\{\begin{array}{ll}
   u_n^{-1}g_i, & n=n_{p(i,2k)},0\leq k<M_i, \cr
   u_n^{-1}{g_i^{-1}}, & n=n_{p(i,2k+1)},0\leq k<M_i,\cr
   u_n^{-1},& \mbox{otherwise}.
  \end{array}\right.
  $$

Then for $q\in\omega$, the following equality holds:
$$(x\oplus w,y\oplus w)|_0^q=\left\{\begin{array}{ll}
 g'_ig_i^{-1}, &  n_{p(i,2k)}\leq q<n_{p(i,2k+1)},0\leq k<M_i, \cr
 1_G,& \mbox{otherwise}.
  \end{array}\right.$$
By (i), we have
$$\lim_q(x\oplus w,y\oplus w)|_0^q=1_G.$$
So $(x\oplus w) E_*(G)(y\oplus w)$, and thus $\psi(x)E_*(H)\psi(y)$ holds. It follows from Lemma~\ref{Lemm basic Tsi property}(2) that
$$\lim_id_H(\psi(x),\psi(y))|_{[l_{p(i,0)},\, l_{p(i,2M_i-1)+1})}=0.$$
For each $j\in\omega$, let $h_j'=S_{n_j}(u_{n_j}x(n_j))$ and $h_j=S_{n_j}(u_{n_j}y(n_j))$.
Note that $S_{n_j}(1_G)=1_H$ for all $j>0$, so
$$\begin{array}{ll}&(\psi(x),\psi(y))|_{l_{p(i,0)}}^{l_{p(i,2M_i-1)+1}-1}\cr
=&h'_{p(i,0)}h'_{p(i,1)}\cdots h'_{p(i,2M_i-1)}(h_{p(i,0)}h_{p(i,1)}\cdots h_{p(i,2M_i-1)})^{-1}.\end{array}$$
Now we have
$$\lim_id_H(h'_{p(i,0)}h'_{p(i,1)}\cdots h'_{p(i,2M_i-1)},h_{p(i,0)}h_{p(i,1)}\cdots h_{p(i,2M_i-1)})=0,$$
and hence $$\lim_i h'_{p(i,0)}h'_{p(i,1)}\cdots h'_{p(i,2M_i-1)}(h_{p(i,0)}h_{p(i,1)}\cdots h_{p(i,2M_i-1)})^{-1}=1_H.$$
Let $V_1\ni 1_{H_{m_0}}$ be an open subset of $H_{m_0}$ with $V_1^2\subseteq V$.
Since $\pi_{m_0}$ is a continuous homomorphism, there exists an $i_0\in\omega$ such that
$$\pi_{m_0}(h'_{p(i,0)}h'_{p(i,1)}\cdots h'_{p(i,2M_i-1)})\pi_{m_0}(h_{p(i,0)}h_{p(i,1)}\cdots h_{p(i,2M_i-1)})^{-1}\in V_1$$
holds for any $i>i_0$.
By $(4)$, we have
$$\lim_i d_{H_{m_0}}(1_{H_{m_0}},\pi_{m_0}(h_{p(i,0)}h_{p(i,1)}\cdots h_{p(i,2M_i-1)}))\leq\lim_i 2^{-i}=0.$$
Therefore, for $i$ large enough,
$$\pi_{m_0}(h_{p(i,0)}h_{p(i,1)}\cdots h_{p(i,2M_i-1)})\in V_1.$$
It follows that
$$\pi_{m_0}(h'_{p(i,0)}h'_{p(i,1)}\cdots h'_{p(i,2M_i-1)})\in V,$$
i.e.,
$$\pi_{m_0}(S_{n_{p(i,0)}}(g'_i)S_{n_{p(i,1)}}((g'_i)^{-1})\cdots S_{n_{p(i,2M_i-2)}}(g'_i)S_{n_{p(i,2M_i-1)}}((g'_i)^{-1}))\in V,$$
contradicting (3).
\hfill\mbox{Claim 1.} $\Box$
\vskip 3mm

For any $m\in\omega$, let $V_m=\{g\in G:d_G(1_G,g)<r_m\}$ and $W_m=\bigcup_iV_m^i$. Note that $V_m^{-1}=V_m$. It is clear that $W_m$ is an open, and thus a clopen subgroup of $G$.
We claim that $W_m\cap \bigcup_n F_n\subseteq X_m$. For $i=0$, note that $1_G\in X_m$, so Claim~1 gives that $V_m\cap \bigcup_n F_n\subseteq X_m$. Assume that $V_m^i\cap \bigcup_n F_n\subseteq X_m$. Let $h\in V_m^{i+1}\cap \bigcup_n F_n$ with $h=gg'$, where $g\in V_m^i$ and $g'\in V_m$. Note that $\bigcup_n F_n$ is a dense subgroup of $G$. We can find $\hat g,\hat g'\in \bigcup_n F_n$ such that $d_G(h,\hat g\hat g')<r_m$ with $\hat g\in V_m^i$ and $\hat g'\in V_m$. Since $d_G(\hat g,\hat g\hat g')=d_G(1_G,\hat g')<r_m$ and $\hat g\in X_m$, we have $\hat g\hat g'\in X_m$, and thus $h\in X_m$. So $W_m\cap \bigcup_n F_n\subseteq X_m$.

Let $g\in W_m\cap \bigcup_n F_n\subseteq X_m$. Pick a $j_g\in\omega$ with $g\in F_{n_{j_g}-1}$. For any $\delta>0$, by the definition of $X_m$, there exists a $j_\delta\geq j_g$ such that
$$\forall j'> j\geq j_{\delta}\,(d_{H_m}(1_{H_m},\pi_m(S_{n_j}(g)S_{n_{j'}}(g^{-1})))<\delta).$$
Let $j,k\in\omega$ with $j_\delta<j<k$. Fix a $k'>k$. Since $\pi_{m}$ is a homomorphism, we have
$$d_{H_m}(1_{H_m},\pi_m(S_{n_j}(g))\pi_m(S_{n_{k'}}(g^{-1})))<\delta,$$
$$d_{H_m}(1_{H_m},\pi_m(S_{n_k}(g))\pi_m(S_{n_{k'}}(g^{-1})))<\delta.$$
So $d_{H_m}(\pi_m(S_{n_j}(g)),\pi_m(S_{n_k}(g)))<2\delta.$ Thus $(\pi_m(S_{n_{j}}(g)))$ is a Cauchy sequence in $H_m$.
By the completeness of $d_{H_m}$, we can define
$$\pi^S_m(g)=\lim_{j}\pi_m(S_{n_{j}}(g))\in H_{m}.$$

Note that for $g,g'\in F_n$ with $n<n_j$,
$$d_H(S_{n_j}(g),S_{n_j}(g'))\le d_H^\infty(T_{n_j}(u_{n_j}^{-1}g),T_{n_j}(u_{n_j}^{-1}g')),$$
and $$\forall h,h'\in H\, (d_{H_m}(\pi_m(h),\pi_m(h'))\leq 2^md_H(h,h')).$$
By Lemma~\ref{continuous property}, $\pi^S_m$ is uniformly continuous on $W_m\cap \bigcup_n F_n$, which can be uniquely extended to a uniformly continuous map from $W_m$ to $H_m$, still denoted by $\pi^S_m$. Let $W=\bigcap_m W_m$, we define a continuous map $S:W\rightarrow H$ as: for $g\in W$,
$$S(g)=(\pi^S_0(g),\pi^S_1(g),\ldots,\pi^S_m(g),\ldots).$$
It is worth noting that $S(1_G)=1_H$ and $W$ is a closed subgroup of $G$.

Now we prove clause (ii) first:

\begin{claim}\label{preserving converges on U}
The map $S^\#$ is a reduction of $E_*(W)$ to $E_*(H)\times E(W;0)$.
\end{claim}

\noindent\emph{Proof of Claim 2.}
Let $x,y\in W^\omega$. Note that $xE_*(G)y\iff xE_*(W)y$ as $W$ is a closed subgroup of $G$. By Proposition~\ref{equiv for sharp}, we only need to show that
$$x E_*(G)y\iff S^\omega(x)E_*(H)S^\omega(y)$$
holds whenever $\lim_p d_G(x(p),y(p))=0$.
Since $\bigcup_n F_n$ is dense, by the definitions of $d_H$ and $S$,
for each $p\in\omega$, we can find a sufficiently large $j(p)\in\omega$ and $g_p,g'_p\in\bigcup_n F_n $ so that

 \begin{enumerate}
   \item [(1)] $0<j(0)<j(1)<\cdots j(p)<\cdots$;
   \item [(2)]$g_p,g'_p\in F_{n_{j(p)}-1}\cap\bigcap_{m\le p} W_m$;
   \item [(3)]$d_G(x(p),g'_p)<2^{-p}$ and $d_G(y(p),g_p)<2^{-p}$; and
   \item [(4)] for all $m\le p$,
   $$d_{H_m}(\pi_m(S(x(p))),\pi_m(S_{n_{j(p)}}(g'_p)))<2^{-(p+1)},$$
   $$d_{H_m}(\pi_m(S(y(p))),\pi_m(S_{n_{j(p)}}(g_p)))<2^{-(p+1)}.$$
 \end{enumerate}

Now we define, for each $n\in I$,
$$\hat x(n)=\left\{\begin{array}{ll}
                         u_n^{-1}g'_p,  & n=n_{j(p)},\cr
                         u_n^{-1}, & \mbox{otherwise,}
                        \end{array}\right. \quad
\hat y(n)=\left\{\begin{array}{ll}
                         u_n^{-1}g_p,  & n=n_{j(p)},\cr
                         u_n^{-1}, & \mbox{otherwise.}
                        \end{array}\right.$$
For any $i<k$, by (3) and Lemma~\ref{Lemm basic Tsi property}(1),
$$d_G((x,y)|_i^k,g'_{i}\cdots g'_{k}{g_k^{-1}}\cdots {g_i^{-1}})<\sum_{q=i}^{k}2^{-q+1}<2^{-i+2}.$$
Then by Lemma~\ref{Lemm basic Tsi property}(2),
$$xE_*(G)y\iff\lim_pg'_{0}\cdots g'_p{g_p^{-1}}\cdots {g_0^{-1}}\mbox{  converges}. $$
For any $k>n_{j(0)}$, let $p_k$ be the largest $p$ with $n_{j(p)}\leq k$. Note that
$$(\hat x\oplus w,\hat y\oplus w)|_0^k=g'_{0}\cdots g'_{p_k}{g_{p_k}^{-1}}\cdots {g_0^{-1}}.$$
It follow that
$$x E_*(G)y\iff(\hat x\oplus w) E_*(G)(\hat y\oplus w)\iff\psi(\hat x)E_*(H)\psi(\hat y).$$

By (4), we have
$$\begin{array}{ll}d_H(S(x(p)),S_{n_{j(p)}}(g'_p))&=\sum_{m=0}^{\infty}2^{-m}d_{H_m}(\pi_m(S(x(p))),\pi_m(S_{n_{j(p)}}(g'_p)))\cr
&<\sum_{m=0}^{p}2^{-(m+p+1)}+\sum_{m=p+1}^{\infty}2^{-m}\cr
&<2^{-(p-1)},\end{array}$$
$$\begin{array}{ll}d_H(S(y(p)),S_{n_{j(p)}}(g_p))&=\sum_{m=0}^{\infty}2^{-m}d_{H_m}(\pi_m(S(y(p))),\pi_m(S_{n_{j(p)}}(g_p)))\cr
&<\sum_{m=0}^{p}2^{-(m+p+1)}+\sum_{m=p+1}^{\infty}2^{-m}\cr
&<2^{-(p-1)}.\end{array}$$
Note that $S_{n_j}(1_G)=1_H$ for all $j\in\omega$, by similar arguments as above, we see that
$$\begin{array}{ll}& S^\omega(x)E_*(H)S^\omega(y)\cr
\iff & \lim_p S_{n_{j(0)}}(g'_0)\cdots S_{n_{j(p)}}(g'_p)S_{n_{j(p)}}(g_p)^{-1}\cdots S_{n_{j(0)}}(g_0)^{-1}\mbox{ converges.}\end{array}$$
For any $i> j(0)$, let $p_{i}$ be the largest $p$ with $j(p)\leq i$. Then we have
$$(\psi(\hat x),\psi(\hat y))|_0^{l_{i+1}-1}=S_{n_{j(0)}}(g'_0)\cdots S_{n_{j(p_{i})}}(g'_{p_{i}})S_{n_{j(p_{i})}}(g_{p_{i}})^{-1}\cdots S_{n_{j(0)}}(g_0)^{-1}.$$
This implies that
$$\lim_i(\psi(\hat x),\psi(\hat y))|_0^{l_{i+1}-1}\mbox{ converges}\iff S^\omega(x)E_*(H)S^\omega(y).$$
Note that $\lim_pd_G(g'_p,g_p)=0$, it follows from Lemma~\ref{preserving converges on F_n} that
$$\lim_i(\psi(\hat x),\psi(\hat y))|_0^{l_{i+1}-1}\mbox{ converges}\iff \psi(\hat x) E_*(H)\psi(\hat y).$$
So $S^\omega(x)E_*(H)S^\omega(y)\iff\psi(\hat x) E_*(H)\psi(\hat y)\iff xE_*(G)y.$
\hfill\mbox{Claim 2.} $\Box$
\vskip 3mm

Subsequently, we prove clause (i) as follows:

\begin{claim}\label{claim_T_m}
For $m\in\omega$, if $x,y\in W_m^\omega$, then
  $$x E_*(G)y\Longrightarrow(\pi^S_m)^\omega(x)E_*(H_m)(\pi^S_m)^\omega(y).$$
\end{claim}

\noindent\emph{Proof of Claim 3.}
Let $x,y\in W_m^\omega$ with $x E_*(G)y$. Then for each $p\in\omega$, we can find a $j(p)\in\omega$ and $g_p,g_p'\in W_m\cap\bigcup_n F_n$ so that
\begin{enumerate}
   \item [(1)] $0<j(0)<j(1)<\cdots j(p)<\cdots$;
   \item [(2)]$g_p,g'_p\in F_{n_{j(p)}-1}\cap W_m$;
   \item [(3)]$d_G(x(p),g'_p)<2^{-p}$ and $d_G(y(p),g_p)<2^{-p}$; and
   \item [(4)] we have
   $$d_{H_m}(\pi^S_m(x(p)),\pi_m(S_{n_{j(p)}}(g'_p)))<2^{-p},$$
 $$d_{H_m}(\pi^S_m(y(p)),\pi_m(S_{n_{j(p)}}(g_p)))<2^{-p}.$$
\end{enumerate}

Now define, for each $n\in I$,
$$\hat x(n)=\left\{\begin{array}{ll}
                         u_n^{-1}g'_p,  & n=n_{j(p)},\cr
                         u_n^{-1}, & \mbox{otherwise,}
                        \end{array}\right. \quad
\hat y(n)=\left\{\begin{array}{ll}
                         u_n^{-1}g_p,  & n=n_{j(p)},\cr
                         u_n^{-1}, & \mbox{otherwise.}
                        \end{array}\right.$$
Following a similar approach to the proof of Claim~\ref{preserving converges on U}, we obtain
$$ x E_*(G)y\iff(\hat x\oplus w) E_*(G)(\hat y\oplus w)\iff\psi(\hat x)E_*(H)\psi(\hat y).$$
So we have
$$\lim_p(\psi(\hat x),\psi(\hat y))|_0^{l_{j(p)+1}-1}\mbox{ converges},$$
and thus
$$\lim_p S_{n_{j(0)}}(g'_0)\cdots S_{n_{j(p)}}(g'_{p})S_{n_{j(p)}}(g_{p})^{-1}\cdots S_{n_{j(0)}}(g_0)^{-1}\mbox{ converges}.$$
Since $\pi_m$ is a continuous homomorphism, we obtain
$$\lim_p \pi_m(S_{n_{j(0)}}(g'_0)\cdots S_{n_{j(p)}}(g'_{p})S_{n_{j(p)}}(g_{p})^{-1}\cdots S_{n_{j(0)}}(g_0)^{-1})\mbox{ converges}.$$
Following again the similar arguments as in the proof of Claim~\ref{preserving converges on U}, we get $(\pi^S_m)^\omega(x)E_*(H_m)(\pi^S_m)^\omega(y)$.
\hfill\mbox{Claim 3.} $\Box$
\vskip 3mm

It follows from Claim~\ref{preserving converges on U} and Claim~\ref{claim_T_m} that $W_m,\pi^S_m,W$, and $S$ are as desired. Finally, by Lemma~\ref{converse}, we can conclude that the converse is also true if $G=W$ and if the interval $[0,1]$ can be embedded into $H$.
\end{proof}

\begin{remark}\label{remark_Pre-rigid Theorem}
Recall that $G_0$ is the connected component of $1_G$ in $G$. Since each $W_m$ in the preceding theorem is an open subgroup of $G$, it is also clopen in $G$. So $G_0\subseteq W_m$ as it is connected, and thus $G_0\subseteq W$.

From the proof of the preceding theorem, it follows that $W$ can be chosen to be clopen when $H$ is TSI strongly NSS. Similarly, by~\cite[Theorem~6.13]{DZ}, $W$ can be also clopen when $H$ is locally compact, and $W$ can be chosen to be $W=G$ when $H$  is compact.
\end{remark}

\section{Rigid theorems}

In this section, we use the Pre-rigid Theorem to prove several Rigid theorems for various classes of TSI Polish groups.

\begin{lemma}\label{homomorphism}
Let $G,H$ be two Polish groups and $S:G\rightarrow H$ a continuous map with $S(1_G)=1_H$. Suppose $H$ is NSS. Then the following are equivalent:

\begin{enumerate}
  \item There exists an open subgroup $W$ of $G$ such that the map $S\upharpoonright W:W\rightarrow H$  is a continuous homomorphism.
  \item There exists an open subgroup $W'$ of $G$ such that the map $S^\omega$ is an $(E_*(W'),E_*(H))$-homomorphism.
\end{enumerate}
\end{lemma}

\begin{proof}
(1)$\Rightarrow$(2). Let $W$ be an open subgroup of $G$ such that $S\upharpoonright W:W\rightarrow H$ is a continuous homomorphism. Then
$$\forall x,y\in W^\omega\,\forall g\in W\,(\lim_n(x,y)|_0^n=g\Rightarrow \lim_n(S^\omega(x),S^\omega(y))|_0^n=S(g)).$$
Since $W$ is also closed, we see that $S^\omega$ is an $(E_*(W),E_*(H))$-homomorphism. Hence $W'=W$ is as required.

(2)$\Rightarrow$(1). Let $W'$ be an open subgroup of $G$ such that the map $S^\omega$ is an $(E_*(W'),E_*(H))$-homomorphism,
and let $d_G,d_H$ be compatible left-invariant metrics on $G$ and $H$ respectively.

Define $X=\{g\in W':S(g^{-1})=S(g)^{-1}\}$. Then $X\ne\emptyset$, since $1_G\in X$.

\setcounter{claim}{0}
\renewcommand{\theclaim}{\arabic{claim}}

\begin{claim}\label{claim g g- for S}
There exists an $r_0>0$ such that
$$\forall g,g'\in W'\, ((d_G(g,g')<r_0\wedge g\in X)\Rightarrow g'\in X).$$
\end{claim}

\noindent\emph{Proof of Claim 1.}
If not, then we can find two sequences $(g_q)$ and $(g'_q)$ in $W'$ so that
\begin{enumerate}
  \item [(\romannumeral 1)] $g_q\in X$ and $\,g'_q\notin X$;
  \item [(\romannumeral 2)] $\lim_q (g'_q)^{-1}g_q=1_G.$
\end{enumerate}
Let $h_q=S((g'_q)^{-1})S(g'_q)$, then $h_q\neq 1_H$ for $q\in\omega$. Since $H$ has no small subgroups, there exists some $D>0$ such that, for each $q\in\omega$, we can find an $m_q>0$ with $ d_{H}(h_q^{m_q},1_{H})\geq D$.

Put $M_{-1}=0$ and $M_q=m_0+\cdots+m_q$ for $q\in\omega$. For each $p\in\omega$, define
$$x(p)=\left\{\begin{array}{ll}(g_q')^{-1}, & p=2(M_{q-1}+i),0\leq i<m_q,\cr g'_q, & p=2(M_{q-1}+i)+1,0\leq i<m_q,\end{array}\right.$$
$$y(p)=\left\{\begin{array}{ll}g_q^{-1}, & p=2(M_{q-1}+i),0\leq i<m_q,\cr g_q, & p=2(M_{q-1}+i)+1,0\leq i<m_q.\end{array}\right.$$
For any $k\in\omega$, let $q_k$ be the largest $q$ such that $M_{q-1}\leq k$. Clearly,
$$(x,y)|_0^{2k+1}=1_G,\quad (x,y)|_0^{2k}=(g'_{q_k})^{-1}g_{q_k}.$$
By (ii), we have $x E_*(G)y$. Since $S^\omega$ is an $(E_*(G),E_*(H))$-homomorphism, $S^\omega(x)E_*(H)S^\omega(y)$ holds, i.e.,
$$\lim_k(S^\omega(x),S^\omega(y))|_0^k\mbox{ converges}.$$
Then by Lemma~\ref{Lemm basic Tsi property}(2) and $g_q\in X$, we have
$$\lim_qh_q^{m_q}=\lim_q(S^\omega(x),S^\omega(y))|_{2M_{q-1}}^{2M_q-1}=1_H.$$
A contradiction!
\hfill\mbox{Claim 1.} $\Box$
\vskip 3mm

Put $V=\{g\in W':d_G(1_G, g)<r_0\}$. Note that $V=V^{-1}$. Define $W_1=\bigcup_iV^i$, which is an open subgroup of $W'$. We claim that $W_1\subseteq X$. Note that $1_G\in X$, so by Claim~\ref{claim g g- for S}, we have $V\subseteq X$. Assume that $V^i\subseteq X$. For any $g=g_1g_2\in V^{i+1}$ with $g_1\in V^i$ and $g_2\in V$, we note that $g_1\in X$ and $d_G(g_1,g)=d_G(g_1,g_1g_2)=d_G(1_G,g_2)<r_0$. Again by Claim~\ref{claim g g- for S}, we have $g\in X$. This shows that $W_1\subseteq X$.

\begin{claim}\label{claim_g'gg-g'-}
There exists $0<r_1<r_0$ such that
$$\forall g,g'\in W_1\,(d_G(1_G,g)<r_1\Rightarrow S(gg')=S(g)S(g')).$$
\end{claim}

\noindent\emph{Proof of Claim 2.}
If not, there are two sequences $(g_q),(g'_q)$ in $W_1$ such that $\lim_qg_q=1_G$ and $S(g_qg'_q)\ne S(g_q)S(g_q')$ for $q\in\omega$.

For each $q\in\omega$, let $h_q=S(g_qg'_q)S((g'_q)^{-1})S(g_q^{-1})$. It follows from $W_1\subseteq X$ that
$$h_q=S(g_qg_q')S(g_q')^{-1}S(g_q)^{-1}\ne 1_H.$$
Since $H$ has no small subgroups, there exists some $D>0$ such that, for each $q\in\omega$, we can find an $m_q\in\omega$ with $ d_{H}(h_q^{m_q},1_{H})\geq D$.

Put $M_{-1}=0$ and $M_q=m_0+\cdots+m_q$ for $q\in\omega$. For each $p\in\omega$, define
$$x(p)=\left\{\begin{array}{ll}g_qg_q', & p=3(M_{q-1}+i),0\leq i<m_q,\cr (g'_q)^{-1}, & p=3(M_{q-1}+i)+1,0\leq i<m_q, \cr g_q^{-1}, &p=3(M_{q-1}+i)+2,0\leq i<m_q,\end{array}\right.$$
$$y(p)=\left\{\begin{array}{ll}g'_q, & p=3(M_{q-1}+i), 0\leq  i<m_q,\cr (g'_q)^{-1}, & p=3(M_{q-1}+i)+1,0\leq i<m_q, \cr 1_G, &p=3(M_{q-1}+i)+2,0\leq i<m_q.\end{array}\right.$$
For any $k\in\omega$, let $q_k$ be the largest $q$ such that $M_{q-1}\leq k$. Clearly,
$$(x,y)|_0^{3k}=(x,y)|_0^{3k+1}= g_{q_k},\quad (x,y)|_0^{3k+2}=1_G.$$
Note that $\lim_qg_q=1_G$, so $x E_*(G)y$ holds. Since $S^\omega$ is an $(E_*(G),E_*(H))$-homomorphism, $S^\omega(x)E_*(H)S^\omega(y)$ holds, i.e.,
$$\lim_k(S^\omega(x),S^\omega(y))|_0^k\mbox{ converges}.$$
By Lemma~\ref{Lemm basic Tsi property}(2) and $g_q'\in W_1\subseteq X$, we have $S(g_q')S((g_q')^{-1})=1_H$. So
$$\lim_qh_q^{m_q}=\lim_q(S^\omega(x),S^\omega(y))|_{3M_{q-1}}^{3M_q-1}=1_H.$$
A contradiction!
\hfill\mbox{Claim 2.} $\Box$
\vskip 3mm

Finally, let $V_1=\{g\in W': d_G(1_G,g)<r_1\}$ and $W=\bigcup_iV^i_1$. Note that $V_1=V_1^{-1}$ and $V_1\subseteq V$, so $W$ is an open subgroup of $W_1$. For any $g,g'\in W$, we shall check that $S(gg')=S(g)S(g')$. There are $g_0,g_1,\cdots,g_m\in V_1$ such that $g=g_0 g_1\cdots g_m$. Note that $W\subseteq W_1$, so by Claim~\ref{claim_g'gg-g'-}, we have
$$S(gg')=S(g_0)S(g_1\cdots g_mg')=S(g_0)\cdots S(g_m)S(g')=S(g)S(g').$$
Then $W$ is as required.
\end{proof}

Consider a Polish group $G$ and a sequence $(g_n)$ in $G$. Recall that $(g_n)$ is \textit{$\iota$-Cauchy} if it is $d$-cauchy for some compatible left-invariant metric $d$ on $G$. This definition is independent of the choice of $d$ (cf.~\cite[Proposition 3.B.1]{becker}).

Let $G,H$ be two Polish groups and $\varphi:G\rightarrow H$ a continuous homomorphism. We define ${\rm IPC(\varphi)}$ as the set of all $x\in G^\omega$ that satisfies:
$$(x(0)\cdots x(p))\mbox{ is $\iota$-Cauchy}\iff(\varphi(x(0))\cdots \varphi(x(p)))\mbox{ is $\iota$-Cauchy}.$$

\begin{lemma}\label{map ker}
Let $G,H$ be two Polish groups and $\varphi:G\rightarrow H$ a continuous homomorphism.
If $x\in {\rm IPC(\varphi)}$ for any $x\in G^\omega$ with $\lim_px(p)=1_G$, then $\ker(\varphi)$ is non-archimedean.

In particular, if $\varphi(G)$ is closed in $H$, then the converse is also true.
\end{lemma}

\begin{proof}
Let $d_G,d_H$ be compatible left-invariant metrics on $G$ and $H$ respectively. Define $V_k=\{g\in\ker(\varphi):d_G(1_G,g)<2^{-k}\}$.

First, suppose $x\in {\rm IPC(\varphi)}$ for any $x\in G^\omega$ with $\lim_px(p)=1_G$.

Assume toward a contradiction that $\ker(\varphi)$ is not non-archimedean. Then there exists a $K\in\omega$ such that $V_{K}$ contains no open subgroups of $\ker(\varphi)$. Since $\bigcup_m V_k^m$ is an open subgroup of $\ker(\varphi)$, we have $\bigcup_n V_k^m\nsubseteq V_{K}$. Thus for each $k\in\omega$, we can find an $m_k\in\omega$ and $g_{k,0},\dots,g_{k,m_k-1}\in V_k$ such that
$g_{k,0}g_{k,1}\cdots g_{k,m_k-1}\notin V_{K}$.

Put $M_{-1}=0$ and $M_k=m_0+\dots+m_k$ for $k\in\omega$. Let $x\in G^\omega$ be defined as $x(p)=g_{k,i}$ for $p=M_{k-1}+i$ with $0\leq i<m_k$. Note that $\lim_px(p)=1_G$, so $x\in {\rm IPC(\varphi)}$. Since $d_G$ is left-invariant, we have
$$\forall k\, (d_G(x(0)\cdots x(M_{k-1}-1),x(0)\cdots x(M_k-1))\geq 2^{-K}).$$
Thus $(x(0)\cdots x(p))$ is not $\iota$-Cauchy.
But $(\varphi(x(0))\cdots\varphi(x(p)))$ is $\iota$-Cauchy since $\varphi(x(p))=1_H$ for all $p\in\omega$, contradicting that $x\in{\rm IPC}(\varphi)$.

Now assume that $\varphi(G)$ is closed in $H$ and $\ker(\varphi)$ is non-archimedean. Let $x\in G^\omega$ with $\lim_px(p)=1_G$. We prove $x\in {\rm IPC(\varphi)}$ as follows:

On the one hand, suppose that $(x(0)\cdots x(p))$ is $\iota$-Cauchy. For any $\varepsilon>0$, since $\varphi$ is continuous, there is a $\delta>0$ such that
$$\forall g\in G\,(d_G(1_G,g)<\delta\Rightarrow d_H(1_H,\varphi(g))<\varepsilon).$$
Since $(x(0)\cdots x(p))$ is $\iota$-Cauchy, we can find an $N\in\omega$ such that
$$\forall m\geq n>N\, (d_G(1_G,x(n)\cdots x(m))<\delta).$$
Thus we have that, for any $ m\geq n>N$,
$$d_H(1_H,\varphi(x(n))\cdots\varphi(x(m)))=d_H(1_H,\varphi(x(n)\cdots x(m)))<\varepsilon.$$
This shows that $(\varphi(x(0))\cdots\varphi(x(p)))$ is $\iota$-Cauchy.

On the other hand, suppose that $(\varphi(x(0))\cdots\varphi(x(p)))$ is $\iota$-Cauchy. Let $(W_p)$ be a decreasing neighborhood basis of $1_{\ker(\varphi)}$ such that each $W_p$ is an open subgroup of $\ker(\varphi)$. Let $\widetilde{\varphi}:G/{\ker(\varphi)}\rightarrow H$ be defined as
$$\widetilde{\varphi}(\ker(\varphi)g)=\varphi(g).$$
Note that $\varphi(G)$ is closed, so it is also a Polish group under the topology inherited from $H$.
Thus $\widetilde{\varphi}$ is a topological isomorphism of $G/{\ker(\varphi)}$ onto $\varphi(G)$ (cf.~\cite[Corollary 2.3.4]{gaobook_1}).
Note that $\phi=\widetilde{\varphi}^{-1}\circ \varphi$ is the canonical projection map, where $\phi(g)=\ker(\varphi)g$ for $g\in G$.
Since $\widetilde{\varphi}$ is a topological isomorphism, we have $(\phi(x(0))\cdots\phi(x(p)))$ is $\iota$-Cauchy.

Let $d_{\phi}(\ker(\varphi)g,\ker(\varphi)g')=\inf\{d_G(hg,h'g'):h,h'\in\ker(\varphi)\}$, then $d_{\phi}$ is a compatible left-invariant metric on $G/\ker(\varphi)$ (cf.~\cite[Lemma 2.2.8]{gaobook_1}).

For any $\varepsilon>0$, we can find $0<\varepsilon'<\varepsilon$ and $p_0\in\omega$ satisfying
$$\{g\in\ker(\varphi):d_G(1_G,g)<3\varepsilon'\}\subseteq W_{p_0}\subseteq\{g\in\ker(\varphi):d_G(1_G,g)<\varepsilon\}.$$
Note that $\lim_px(p)=1_G$ and $(\phi(x(0))\cdots\phi(x(p)))$ is $\iota$-Cauchy. Thus there is an $N\in\omega$ such that
$$d_G(1_G,x(n))<\varepsilon',$$
$$d_{\phi}(\phi(x(N)\cdots x(n)),\ker(\varphi)1_H)<\varepsilon'$$
for each $n>N$. By the definition of $d_\phi$, there exists some $g_n\in\ker(\varphi)$ for each $n>N$ such that
$$d_G(x(N)\cdots x(n),g_n)<\varepsilon'.$$
For any $n>N$, we have
$$d_G(x(N)\cdots x(n),x(N)\cdots x(n+1))=d_G(1_G,x(n+1))<\varepsilon'.$$
All there together imply that $d_G(1_G,g_n^{-1}g_{n+1})=d_G(g_n,g_{n+1})<3\varepsilon'$. So $g_n^{-1}g_{n+1}\in W_{p_0}$. Since $W_{p_0}$ is a subgroup, it follows that $g_n^{-1}g_m\in W_{p_0}$ for any $m\geq n\geq N$. This gives that
$$\begin{aligned}
d_G(x(0)\cdots x(n),x(0)\cdots x(m))&=d_G(x(N)\cdots x(n),x(N)\cdots x(m))\cr
&<2\varepsilon'+d_G(g_n,g_m)\cr
&<2\varepsilon'+\varepsilon<3\varepsilon.
\end{aligned}
$$
Therefore, $(x(0)\cdots x(p))$ is $\iota$-Cauchy.
\end{proof}

\begin{proposition}\label{IPC iff reduction}
Let $G,H$ be two TSI Polish groups and $\varphi:G\to H$ a continuous homomorphism. Then the following are equivalent:
\begin{enumerate}
\item[(1)] for $x\in G^\omega$, if $\lim_px(p)=1_G$, then $x\in{\rm IPC}(\varphi)$; and
\item[(2)] the map $\varphi^\#$ is a reduction of $E_*(G)$ to $E_*(H)\times E(G;0)$.
\end{enumerate}
\end{proposition}

\begin{proof}
Since $G$ and $H$ are both TSI, for $x\in {G}^\omega$, we have
$$(x(0)\cdots x(p))\mbox{ is $\iota$-Cauchy}\iff xE_*(G)1_{G^\omega},$$
$$(\varphi(x(0))\cdots\varphi(x(p)))\mbox{ is $\iota$-Cauchy}\iff \lim_p\varphi^\omega(x)E_*(H)1_{H^\omega}.$$
So (2)$\Rightarrow$(1) follows from Proposition~\ref{equiv for sharp}. We prove (1)$\Rightarrow$(2) as follows:

Let $d_G$ be a compatible compete two-sided invariant metric on $G$.

Given $x,y\in G^\omega$. Suppose $xE_*(G)y$. Since $\varphi$ is a continuous homomorphism, we have $\varphi^\omega(x)E_*(H)\varphi^\omega(y)$ holds and $\lim_pd_G(x(p),y(p))=0$.

On the other hand, suppose $\varphi^\omega(x)E_*(H)\varphi^\omega(y)$ holds and $\lim_pd_G(x(p),y(p))=0$. Define $z\in G^\omega$ as:
$$z(p)=y(0)\cdots y(p-1)x(p)y(p)^{-1}y(p-1)^{-1}\cdots y(0)^{-1}$$
for $p\in\omega$. Since $d_G$ is two-sided invariant,
$$\lim_pd_G(z(p),1_G)=\lim_pd_G(x(p),y(p))=0.$$
So $z\in{\rm IPC}(\varphi)$. Since $\varphi$ is homomorphism, it holds that
$$\varphi(z(p))=\varphi(y(0))\cdots\varphi(y(p-1))\varphi(x(p))\varphi(y(p))^{-1}\varphi(y(p-1))^{-1}\cdots\varphi(y(0))^{-1}.$$
Therefore,
$$xE_*(G)y\iff(z(0)\cdots z(p))\mbox{ is $\iota$-Cauchy},$$
$$\varphi^\omega(x)E_*(H)\varphi^\omega(y)\iff(\varphi(z(0))\cdots\varphi(z(p)))\mbox{ is $\iota$-Cauchy}.$$
So (1)$\Rightarrow$(2) follows from Proposition~\ref{equiv for sharp} again.
\end{proof}

Now we are ready to prove the following theorem, which is crucial for the rest of this article.

\begin{theorem}\label{Theore TSI strongly NSS homomorphism}
Let $G$ be a TSI Polish group, and let $H$ be a closed subgroup of the product of a sequence of TSI strongly NSS Polish groups $(H_m)$. If $E(G)\le_BE(H)$, then there exists a continuous homomorphism $S:W\rightarrow H$ such that $\ker(S)$ is non-archimedean, where $W\supseteq G_0$ is a countable intersection of clopen subgroups in $G$.

In particular, the converse holds if $G=W$, $S(G)$ is closed in $H$, and the interval $[0,1]$ can be embedded into $H$.
\end{theorem}

\begin{proof}
Suppose $E(G)\le_BE(H)$. It follows from Theorem~\ref{Pre-rigid Theorem} that, for each $m\in\omega$, there exists an open subgroup $W_m'$ of $G$ and a continuous map $\pi^S_m:W_m'\to H_m$ with $\pi^S_m(1_G)=1_{H_m}$ such that (1) $(\pi^S_m)^\omega:(W_m')^\omega\rightarrow H_m^\omega$ is an $(E_*(W_m'),E_*(H_m))$-homomorphism; and
(2) let $W'=\bigcap_mW_m'$ and define $S:W'\rightarrow H$ as $S(g)=(\pi^S_0(g),\pi^S_1(g),\cdots)$,
then $S^\#$ is a reduction of $E_*(W')$ to $E_*(H)\times E(W';0)$.

Note that $H_m$ is NSS. By Lemma~\ref{homomorphism}, there is an open subgroup $W_m$ of $W_m'$ such that $\pi^S_m\upharpoonright W_m$ is a homomorphism.
We put $W=\bigcap_m W_m$, and still denote $S\upharpoonright W$ by $S$ for brevity.
Then $S$ is a continuous homomorphism. It is clear that $G_0\subseteq W\subseteq W'$.

As both $G$ and $H$ are TSI, according to Lemma~\ref{map ker} and Proposition~\ref{IPC iff reduction}, $\ker(S)$ is non-archimedean. Thus, $W$ and $S$ fulfill the requirements.

On the other hand, let $S$ be a continuous homomorphism from $G$ to $H$ such that $\ker(S)$ is non-archimedean. Suppose that $S(G)$ is closed in $H$ and the interval $[0,1]$ can be embedded into $H$. Again by Lemma~\ref{map ker} and Proposition~\ref{IPC iff reduction}, the map $S^\#$ is a reduction of $E_*(G)$ to $E_*(H)\times E(G;0)$.
Finally, it follows from Lemma~\ref{converse} that $E(G)\leq_B E(H)$.
\end{proof}

\begin{remark}\label{remark nss homor}
Let $G$ be a TSI Polish group and $H$ a TSI strongly NSS Polish group. Suppose that $E(G)\leq_B E(H)$. Then it follows from the proof of Theorem~\ref{Theore TSI strongly NSS homomorphism} that there exist an open subgroup $W\supseteq G_0$ of $G$ and a continuous homomorphism $S:W\rightarrow H$ such that $\ker(S)$ is non-archimedean.
Moreover, the map $S^\#$ is a reduction of $E_*(W)$ to $E_*(H)\times E(W;0)$.

This observation will be crucial in section 7.
\end{remark}

The following is an immediate corollary.

\begin{corollary}\label{cor totall disc}
Let $G$ be a TSI Polish group and $H$ a closed subgroup of the product of a sequence of TSI strongly NSS Polish groups. If $H$ is totally disconnected but $G$ is not, then $E(G)\nleq_BE(H)$.
\end{corollary}

\begin{proof}
Assume for contradiction that $E(G)\leq_B E(H)$. By Theorem~\ref{Theore TSI strongly NSS homomorphism}, there is a continuous homomorphism $S:G_0\rightarrow H$ such that $\ker(S)$ is non-archimedean.
Note that $S(G_0)$ is also connected, so $S(G_0)\subseteq H_0$. Since $H$ is totally disconnected, we have $S(G_0)=\{1_H\}$. So $G_0=\ker(S)$, which is connected and non-archimedean. This implies that $G_0=\{1_G\}$, contradicting that $G$ is not totally disconnected.
\end{proof}

Note that the converses of Lemma~\ref{map ker} and Theorem~\ref{Theore TSI strongly NSS homomorphism} require $S(G)$ to be closed in $H$.
Next we present a lemma in which this requirement can be avoided.

\begin{lemma}\label{locally homeomorphism}
Let $G,H$ be two TSI Polish groups. Suppose $S:G\to H$ is a continuous homomorphism, and $U\ni 1_G$ is an open subset of $G$ such that $S\upharpoonright U:U\rightarrow S(U)$ is a homeomorphism.
Then the map $S^\#$ is a reduction of $E_*(G)$ to $E_*(H)\times E(G;0)$.
\end{lemma}

\begin{proof}
Let $d_G,d_H$ be compatible complete two-sided invariant metrics on $G$ and $H$ respectively. By Proposition~\ref{equiv for sharp}, we only need to show that
$$xE_*(G)y\iff S^\omega(x)E_*(H)S^\omega(y)$$
holds whenever $\lim_pd_G(x(p),y(p))=0$.

Suppose $xE_*(G)y$ holds. Since $S:G\to H$ is a continuous homomorphism, it follows trivially that $S^\omega(x)E_*(H)S^\omega(y)$.

On the other hand, suppose $\lim_pd_G(x(p),y(p))=0$ and $S^\omega(x)E_*(H)S^\omega(y)$. We prove $xE_*(G)y$ as follows:

If not, by Lemma~\ref{Lemm basic Tsi property}(2), $\lim_p\sup_{q>p}{d_G(x,y)|_{[p,q]}}\ne 0$. So there exist $r>0$ and two strictly increasing sequences of natural numbers $(p_k),(q_k)$ such that $p_k<q_k$ and $d_G(1_G,(x,y)|_{p_k}^{q_k})=d_G(x,y)|_{[p_k,q_k]}>r$ for each $k\in\omega$.
Let $0<r_0<r$ satisfying that $\{g\in G:d_G(1_G,g)<r_0\}\subseteq U$. There exists a $K\in\omega$ such that $d_G(x(p),y(p))<r_0/2$ for $p>p_K$.
For $p>p_k$ and $k>K$, since $d_G$ is two-sided invariant, we have
$$d_G((x,y)|_{p_k}^p,(x,y)|_{p_k}^{p-1}))=d_G(x(p)y(p)^{-1},1_G)<r_0/2,$$
and hence
$$d_G(1_G,(x,y)|_{p_k}^p)<d_G(1_G,(x,y)|_{p_k}^{p-1})+r_0/2.$$
For each $k>K$, we can find $p_k<p_k'<q_k$ such that
$$r_0/2\le d_G(1_G,(x,y)|_{p_k}^{p_k'})<r_0,$$
so $(x,y)|_{p_k}^{p_k'}\in U$. By $S^\omega(x)E_*(H)S^\omega(y)$, it is clear that
$$S((x,y)|_{p_k}^{p_k'})=(S^\omega(x),S^\omega(y))|_{p_k}^{p_k'}\to 1_H.$$
Since $S\upharpoonright U$ is a homeomorphism from $U$ to $S(U)$, we have $(x,y)|_{p_k}^{p_k'}\to 1_G$. This contradicts that $d_G(1_G,(x,y)|_{p_k}^{p_k'})\geq r_0/2$.
\end{proof}

\subsection{Applications on Lie groups, locally compact groups, and pro-Lie groups}

Recall that a \textit{Lie group} is a group which is also a smooth manifold such that the group operations are smooth functions.
A topological group is a Lie group iff it is locally compact NSS (cf. Page 159 of~\cite{HMS}).
Clearly, a Lie group is Polish iff it is separable iff it has only countably many connected components.
Let $G$ be a Lie group, then $G_0$ is an open normal subgroup of $G$. For more details on Lie groups, we refer to~\cite{varadarajan}.

A completely metrizable topological group $G$ is called a \textit{pro-Lie group} if every open neighborhood of $1_G$ contains a normal subgroup $N$ such that $G/N$ is a Lie group (cf.~\cite[Definition 1]{HM07}). For more details on pro-Lie groups, we refer to~\cite{HM07}.

Applying Theorem~\ref{Theore TSI strongly NSS homomorphism}, we obtain the following result:

\begin{theorem}\label{local compa tsi}
Let $G,H$ be two TSI Polish groups such that $H$ is a pro-Lie group. If $E(G)\leq_B E(H)$, then there exists a continuous homomorphism $S:G_0\rightarrow H$ such that $\ker(S)$ is non-archimedean.

In particular, the converse holds if $G$ is connected and $S(G)$ is closed in $H$.
\end{theorem}

\begin{proof}
Since $H$ is a pro-Lie group, there exist a countable open neighborhood basis $(U_m)$ of $1_H$ and a sequence of normal subgroups $(N_m)$ of $H$ such that $N_m\subseteq U_m$ and $H/N_m$ is a Lie group for each $m\in\omega$. Note that $N_m$ is closed in $H$, since $H/N_m$ is Hausdorff. Without loss of generality, assume that $N_m\supseteq N_{m+1}$.
We define a map $f:H\rightarrow \prod_mH/N_m$ as $f(h)(m)=hN_m$ for $h\in H,m\in\omega$. By~\cite[Proposition 2.3]{HMS1}, $f$ is a topologically isomorphic embedding and $f(H)$ is closed in $\prod_mH/N_m$.

By~\cite[Exercise 2.2.8]{gaobook_1}, all $H/N_m$ are separable TSI Lie groups. As mentioned above, all Lie groups are locally compact NSS.
By Proposition~\ref{NSS and strongly NSS}(2), all Lie groups are strongly NSS. So the first part of the theorem follows from Theorem~\ref{Theore TSI strongly NSS homomorphism}.

For proving the second part of the theorem, by Theorem~\ref{Theore TSI strongly NSS homomorphism}, we only need to show that the interval $[0,1]$ can be embedded into $H$. Since $S(G)$ is a closed subgroup of $H$, it is also a pro-Lie Polish group. Clearly, $S(G)$ is non-singleton and connected. Note that, for any non-singleton connected pro-Lie Polish group $K$, there exists a nontrivial continuous homomorphism $\gamma:\R\rightarrow K$~(see Proposition 19 and Definition 2.6 of~\cite{HM07}). All of these together allow us to embed the interval $[0,1]$ into $S(G)$, and also into $H$.
\end{proof}

A topological group $G$ is called a \textit{SIN-group} if $G$ admits arbitrarily small invariant identity neighborhoods, or equivalently, $G$ has a neighborhood basis $(U_i)_{i\in I}$ of $1_G$ such that $gU_i g^{-1}=U_i$ for all $g\in G$ (cf.~\cite[Definition 2.1]{HMS}).
It follows from~\cite[Exercise 2.1.4]{gaobook_1} that a Polish group is SIN iff it is TSI.
By ~\cite[Theorem 3.6]{HMS}, every locally compact TSI Polish group is a pro-Lie group.

We point out that, in~\cite[Theorem 1.2]{DZlcoabe}, groups $G$ and $H$ are required to be abelian with $G$ being compact. Following is a generalized theorem that removes these requirements.

\begin{theorem}[Rigid Theorem for locally compact TSI groups]\label{rigid-compact}
Let $G$ be a locally compact connected TSI Polish group, $H$ a TSI pro-Lie Polish group. Then $E(G)\leq_B E(H)$ iff there exists a continuous homomorphism $S:G\rightarrow H$ such that $\ker(S)$ is non-archimedean.
\end{theorem}

\begin{proof}
$(\Rightarrow)$. It follows from Theorem~\ref{local compa tsi} and $G=G_0$.

$(\Leftarrow)$. Let $\varphi:G\to G/\ker(S)$ be the canonical projection. Note that $\varphi$ is a continuous surjective homomorphism with $\varphi(g)=\ker(S)g$ for $g\in G$.
Clearly, $\ker(\varphi)=\ker(S)$ is non-archimedean.
By Lemma~\ref{map ker} and Proposition~\ref{IPC iff reduction}, the map $\varphi^\#$ is a continuous reduction of $E_*(G)$ to $E_*(G/\ker(S))\times E(G;0)$.

Let $S^*:G/\ker(S)\rightarrow H$ be the map defined as $S^*(\ker(S)g)=S(g)$ for $g\in G$. It is clear that $S^*$ is a continuous injective homomorphism.
Let $d_G$ be a compatible complete two-sided invariant metric on $G$, and define
$$d^*(\ker(S)g,\ker(S)g')=\inf\{d_G(kg,k'g'):k,k'\in\ker(S)\}.$$
Then $d^*$ is a compatible metric on $G/\ker(S)$ (cf.\cite[Lemma 2.2.8]{gaobook_1}).
Since $G$ is locally compact, we can find some $r>0$ such that the closure $\overline{V}$ of the open set $V=\{g\in G:d_G(1_G,g)<r\}$ is compact. Let
$$U=\{\ker(S)g:d^*(1_{G/\ker(S)},\ker(S)g)<r\}.$$
Since $d_G$ is two-sided invariant, from the definition of $d^*$, we see that $U\subseteq S(V)$.
Note that $\overline{U}\subseteq S(\overline{V})$, so $\overline{U}$ is also compact.
This implies that $S^*\upharpoonright\overline{U}:\overline{U}\rightarrow S^*(\overline{U})$ is a homeomorphism.
By Lemma~\ref{locally homeomorphism}, the map $(S^*)^\#$ is a continuous reduction of $E_*(G/\ker(S))$ to $E_*(H)\times E(G/\ker(S);0)$.

Given $x,y\in G^\omega$ with $\lim_pd_G(x(p),y(p))=0$. It follows from Proposition~\ref{equiv for sharp} that
$$xE_*(G)y\iff\varphi^\omega(x)E_*(G/\ker(S))\varphi^\omega(y).$$
The continuity of $\varphi$ implies $\lim_pd^*(\varphi(x(p)),\varphi(y(p)))=0$, so we also have
$$\varphi^\omega(x)E_*(G/\ker(S))\varphi^\omega(y)\iff (S^*)^\omega(\varphi^\omega(x))E_*(H)(S^*)^\omega(\varphi^\omega(y)).$$
Note that $S=S^*\circ\varphi$, so $S^\omega=(S^*)^\omega\circ\varphi^\omega$. It is clear that
$$xE_*(G)y\iff S^\omega(x)E_*(H)S^\omega(y).$$
Thus the map $S^\#$ is a continuous reduction of $E_*(G)$ to $E_*(H)\times E(G;0)$.

Since $G$ is connected, $\{1_H\}\ne S(G)\subseteq H_0$. So $H_0$ is non-singleton. As a closed subgroup of $H$, $H_0$ is also a TSI pro-Lie Polish group.
Similar to the last paragraph of the proof of Theorem~\ref{rigid-compact}, we can embed the interval $[0,1]$ into $H_0$, and also in $H$.
Therefore, the $(\Leftarrow)$ part follows from Lemma~\ref{converse}.
\end{proof}

By restricting our analysis to Lie groups, we provide an affirmative response to Question 7.4 of~\cite{DZ} as follows:

\begin{theorem}[Rigid Theorem for TSI Lie groups]
Let $G,H$ be two separable TSI Lie groups such that $G$ is connected. Then $E(G)\leq_B E(H)$ iff there exists a continuous locally injective homomorphism $S:G\rightarrow H$.
\end{theorem}

\begin{proof}
By Theorem~\ref{rigid-compact}, we only need to show that $\ker(S)$ is discrete. Note that any non-archimedean subgroup of a Lie group is discrete, as all Lie groups are NSS.
\end{proof}

\subsection{Applications on Banach spaces and Fr\'echet spaces}

Now we focus on infinite dimensional vector spaces. Let us recall some elementary notions.
All separable \textit{Fr{\'e}chet spaces}, i.e., separable, completely  metrizable topological vector spaces (see~\cite[5-1]{Wil}), can be viewed as abelian Polish groups under the addition operation.
In this article, all vector spaces are assumed to be real, i.e., over the field $\R$. This is because, in accordance with view of Borel reducibility, the equivalence relation induced by a separable Fr{\'e}chet space is independent of the choice of the field of scalars.

Let $X$ be a vector space. A map $\Vert\cdot\Vert:X\rightarrow\R$ is called a \textit{total paranorm} (see~\cite[2-1]{Wil}) if
\begin{enumerate}
  \item $\Vert v\Vert\geq 0$, $\Vert-v\Vert=\Vert v\Vert$, and $\Vert v\Vert=0\iff v=0$;
  \item ${\Vert}v+v'\Vert\leq {\Vert}v\Vert+{\Vert}v'\Vert$;
  \item for $(t_n)\in\R^\omega,(v_n)\in X^\omega$, if $t_n\rightarrow t$ and $\Vert v_n-v\Vert\rightarrow 0$, then $\Vert t_nv_n-tv\Vert\rightarrow 0$.
\end{enumerate}
A total paranorm $\Vert\cdot\Vert$ is called a \textit{norm} if $\Vert tv\Vert=|t|\Vert v\Vert$ for $t\in\R$ and $v\in X$.
Any Fr{\'e}chet space admits a compatible complete two-sided invariant metric $d$, which is given by a total paranorm $\Vert\cdot\Vert$ as $d(v,v')={\Vert}v-v'\Vert$.

In particular, all separable Banach spaces are separable Fr{\'e}chet spaces. The following are some classical separable Banach spaces:
$$c_0,\, C[0,1], \,l_p, \,L_p[0,1]\quad(\forall p\in[1,+\infty)).$$
The classical Banach-Mazur theorem~\cite[II.B]{Woj} asserts that every separable Banach space is isometric to a subspace of $C[0,1]$. Thus for any separable Banach space $X$, we have $E(X)\leq_B E(C[0,1])$.
For more details on Fr{\'e}chet spaces and Banach spaces, we refer to~\cite{Wil,Woj}.

\begin{theorem}[Rigid Theorem for Fr\'echet spaces]\label{Frechet}
Let $X,Y$ be two separable Fr{\'e}chet spaces such that $Y$ is a closed subgroup of the product of a sequence of TSI strongly NSS Polish groups. Then $E(X)\leq_B E(Y)$ iff $X$ is topologically isomorphic to a closed linear subspace of $Y$.
\end{theorem}

\begin{proof}
$(\Leftarrow)$. It follows from Proposition~\ref{pro closed group}.

$(\Rightarrow)$. Suppose $E(X)\leq_B E(Y)$. Note that $X$ is connected. By Theorem~\ref{Pre-rigid Theorem} and Theorem~\ref{Theore TSI strongly NSS homomorphism}, there exists a continuous homomorphism $S:X\to Y$ such that $\ker(S)$ is non-archimedean, and for ${x},{y}\in X^\omega$, if $\lim_p{x}(p)-{y}(p)=0$, then
$$\sum_px(p)-y(p)\mbox{ converges}\iff\sum_pS(x(p))-S(y(p))\mbox{ converges}.$$
For any integers $m,n$ with $n>0$ and $v\in X$, we have that $nS(\frac{1}{n}v)=S(v)$, and thus $S(\frac{mv}{n})=\frac{m}{n}S(v)$. By its continuity, $S$ is a $\R$-linear map.
Note that $\ker(S)$ is a non-archimedean closed linear subspace of $X$, so $\ker(S)=\{0\}$. Therefore $S$ is injective. We will show that $S$ is a topological isomorphism from $X$ onto $S(X)$.

Let $\Vert\cdot\Vert_X$, $\Vert\cdot\Vert_Y$ be total paranorms on $X$ and $Y$ respectively.
Assume for contradiction that there exist a sequence $(v_q)$ in $X$ and a $\delta>0$ such that
\begin{enumerate}
  \item [(1)] $\lim_q \Vert S(v_q)\Vert_Y=0$, and
  \item [(2)]$\inf_q\Vert v_q\Vert_X>\delta$.
\end{enumerate}
By transfering to a subsequence, we may assume that $\Vert S(v_q)\Vert_Y<2^{-q}$ for all $q\in\omega$.
For each $q\in\omega$, we can find an integer $m_q>0$ such that $\Vert\frac{v_q}{m_q} \Vert_X<2^{-q}$.

Put $M_{-1}=0$ and $M_q=m_0+\cdots+m_q$ for $q\in\omega$.
For each $p\in\omega$, define $x\in X^\omega$ with $x(p)=\frac{v_q}{m_q}$ for $p=M_{q-1}+i$ and $0\leq i<m_q$. It follows that $\lim_px(p)=0$, so
$$\sum_px(p)\mbox{ converges}\iff\sum_pS(x(p))\mbox{ converges}.$$

Fix a $k\in\omega$. By~\cite[Exercise E7.12(ii)]{HM13}, there is an open set $V\ni 0$ in $Y$ such that
$$\forall r\in[0,1]\,\forall h\in V\, (\Vert rh\Vert_Y<2^{-(k+2)}).$$
By (1), there is a $Q>k$ such that
$$\forall q\,( q>Q\Rightarrow S(v_q)\in V).$$
For any $p'>p>M_Q$, there are $q'\ge q>Q$ such that $M_{q-1}\leq p< M_{q}$ and $M_{q'-1}\leq p'<M_{q'}$. It follows that
$$S(x(p))+\cdots+ S(x(M_q-1))=\frac{M_{q}-p}{m_q}S(v_q),\quad S(v_q)\in V,$$
$$S(x(M_{q'-1}))+\cdots+ S(x(p'))=\frac{p'-M_{q'-1}+1}{m_{q'}}S(v_{q'}),\quad S(v_{q'})\in V.$$
These imply that
$$\begin{aligned}
\left\Vert\sum_{i=p}^{p'} S(x(i))\right\Vert_Y&\leq \Vert S(x(p))+\cdots+ S(x(M_q-1))\Vert_Y+ \sum_{i=q+1}^{q'-1}\Vert S(v_i)\Vert_Y\cr
&+\Vert S(x(M_{q'-1}))+\cdots+ S(x(p'))\Vert_Y\cr
&\leq 2^{-(k+2)}+\sum_{i=q+1}^{q'-1}2^{-i}+2^{-(k+2)}<2^{-k}.
\end{aligned}
$$
This shows that $(S(x(0))+\cdots+S(x(p)))$ is a Cauchy sequence. So the series $\sum_{p}S(x(p))$ converges.

For any $q\in\omega$, we have
$$\Vert{x}(M_{q-1})+{x}(M_{q-1}+1)+\cdots+{x}(M_{q}-1)\Vert_X=\Vert v_q\Vert_X>\delta,$$
and thus $\sum_{p}{x}(p)\mbox{ diverges}$. A contradiction!

Now we see that
$$\forall (v_q)\in X^\omega\,(\lim_qS(v_q)=0\Rightarrow\lim_qv_q=0).$$
This implies that $S^{-1}:S(X)\rightarrow X$ is continuous. So $S$ is a topological linear isomorphism from $X$ onto $S(X)$, and hence $S(X)$ is closed.
\end{proof}

The following theorem characterizes strongly NSS separable Banach spaces.

\begin{theorem}\label{pro strong nss banach}
A separable Banach space $X$ is not strongly NSS iff it has a closed linear subspace topologically isomorphic to $c_0$.
\end{theorem}

\begin{proof}
$(\Leftarrow)$. It follows from the fact that any closed subgroup of a strongly NSS Polish group is strongly NSS, whereas $c_0$ is not strongly NSS.

$(\Rightarrow)$. Suppose $X$ is not strongly NSS. The norm of $X$ denoted by $\|\cdot\|$. Let $V=\{v\in X: \Vert v\Vert<1\}$. Then $V$ is not an unenclosed set of $X$, and thus there exists a sequence $(v_i)$ in $X$ such that $\lim_i v_i\ne 0$ and
$$\forall j<k\,\forall \theta\in\{-1,1\}^{\{j,\ldots,k\}}\,\Vert \theta(j)v_j+\cdots+\theta(k)v_{k}\Vert<2.$$

For any integer $n>0$ and $z^0,z^1,\ldots,z^m\in\R^n$, let
$${\rm Con}(\{z^0,z^1,\ldots,z^m\})=\left\{\sum_{i=0}^{m}\lambda_iz^i:\sum_{i=0}^{m}\lambda_i=1,\,\forall i\le m\,(\lambda_i\geq 0)\right\}.$$
Define $\mathbb{D}(n)=\{z\in \R^n:\forall i<n\,(|z(i)|\leq 1)\}$. It is clear that
$$\mathbb{D}(n)={\rm Con}(\{\theta\in\R^n:\forall i<n\,(\theta(i)=\pm 1)\}).$$

Now let $(t_i)\in c_0$. We claim that $\sum_{i}t_iv_i$ converges.
For any $\varepsilon>0$, we can find an $i_0\in\omega$ so that $|t_i|<\varepsilon$ for $i>i_0$.
Let $r=\sup\{|t_i|:i>i_0\}$. We may assume $r> 0$.
Then for any $i_0<j<k$, we have $(t_j,\ldots,t_k)/r\in\mathbb{D}(k-j+1)$. So there is a $\lambda:\{-1,1\}^{\{j,\ldots,k\}}\to\R$ with $\sum_{\theta\in\{-1,1\}^{\{j,\ldots,k\}}}\lambda(\theta)=1$ and $\lambda(\theta)\geq 0$ for each $\theta$ such that
$$(t_j,\ldots,t_k)/r=\sum_{\theta\in\{-1,1\}^{\{j,\ldots,k\}}}\lambda(\theta)(\theta(j),\ldots,\theta(k)).$$
So
$$\sum_{i=j}^{k}\frac{t_i}{r}v_i=\sum_{\theta\in\{-1,1\}^{\{j,\ldots,k\}}}\lambda(\theta)\sum_{i=j}^{k}\theta(i)v_i.$$
Then we have
$$\begin{aligned}
\left\Vert\sum_{i=j}^{k}t_iv_i\right\Vert&\leq\sum_{\theta\in\{-1,1\}^{\{j,\ldots,k\}}}r\lambda(\theta)\left\Vert\sum_{i=j}^{k}\theta(i)v_i\right\Vert\cr
&<\sum_{\theta\in\{-1,1\}^{\{j,\ldots,k\}}}2r\lambda(\theta)=2r\leq 2\varepsilon.
  \end{aligned}$$
So by Cauchy Criterion, $\sum_{i}t_iv_i$ converges.

Therefore, $\sum_{i}t_iv_i$ converges whenever $(t_i)\in c_0$. Note that $\sum_{i}v_i$ diverges, since $\lim_i v_i\ne 0$. By Proposition 2.e.4 of~\cite{LT} and its remark, $X$ has a closed linear subspace topologically isomorphic to $c_0$.
\end{proof}

 Applying previous results to separable Banach spaces, we get:

\begin{theorem}[Rigid Theorem for Banach spaces]\label{rigid theo}
Let $X,Y$ be two separable Banach spaces such that $Y$ contains no closed linear subspaces topologically isomorphic to $c_0$. Then $E(X)\leq_B E(Y)$ iff $X$ is topologically isomorphic to a closed linear subspace of $Y$.
\end{theorem}

\begin{proof}
It follows from theorems~\ref{Frechet} and~\ref{pro strong nss banach}.
\end{proof}

\begin{remark}
It is well known that $l_p$ and $L_p[0,1]$ contain no closed linear subspaces topologically isomorphic to $c_0$, where $p\in[1,+\infty)$. Note that $l_p$ is topologically isomorphic to a closed linear subspace of $l_q$ iff $p=q$ (cf.~\cite[Theorem 5.1]{WJS}). Let $X$ be a separable Fr{\'e}chet space and let $Y$ be the space $l_p$ or $L_p[0,1]$. Then $E(X)\leq_B E(Y)$ iff $X$ is topologically isomorphic to a closed linear subspace of $Y$. In particular,
$$E(l_p)\leq E(l_q)\iff p=q.$$
\end{remark}

Without the assuming of strongly NSS for TSI Polish groups $G$ and $H$, we do not know how to compare $E(G)$ and $E(H)$ with respect to Borel reducibility so far.
For instance, it is not known whether $E(C[0,1])\leq_B E(c_0)$, or whether $E(l_p)\leq_BE(c_0)$ for $p\in(0,+\infty)$.

A worthwhile question is:

\begin{question}
Let $X,Y$ be two separable Fr{\'e}chet spaces. Does it hold that $E(X)\leq_B E(Y)$ iff $X$ is topologically isomorphic to a closed linear subspace of $Y$?
\end{question}

\section{Uniformly NSS Polish groups}

Recall that a topological group $G$ is \emph{uniformly NSS} (cf.~\cite[Definition 11]{CR}) if there is an open subset $V\ni 1_G$ of $G$ such that, for any open subset $U\ni 1_G$ of $G$,
$$\exists n\,\forall g\in G\,(g\notin U\Rightarrow\exists m\leq n\,(g^m\notin V)).$$
Clearly, a locally compact Polish group is uniformly NSS iff it is NSS.

A topological group is a \textit{Banach-Lie group} if it is a Banach manifold such that the group operations are smooth functions (cf.~\cite[Section 6]{Upmeier}).
All Banach spaces and Lie groups are Banach-Lie groups, while all Banach-Lie groups are uniformly NSS (cf.~\cite[Theorem 2.7]{Morris}).

Under the assumption of uniformly NSS, we show that every Borel reduction of $E(G)$ to $E(H)$ results in a continuous homomorphism, which is also a locally homeomorphism.

\begin{theorem}\label{unform nss}
Let $G,H$ be two TSI Polish groups. Assume that $G$ is uniformly NSS and $H$ is strongly NSS. If $E(G)\leq_B E(H)$, then there exist an open subgroup $W$ of $G$, an open neighborhood $U\subseteq W$ of $1_G$, and a continuous homomorphism $S:W\to H$ such that $S\upharpoonright\overline{U}:\overline{U}\rightarrow S(\overline{U})$ is a homeomorphism and $S(\overline{U})$ is closed in $H$.

In particular, the converse is also true if $G=W$ and the interval $[0,1]$ can be embedded into $H$.
\end{theorem}

\begin{proof}
Let $d_G,d_H$ be compatible complete two-sided invariant metrics on $G$ and $H$ respectively.
Assume that $E(G)\leq_B E(H)$. It follows from Remark~\ref{remark nss homor} that there exist an open subgroup $W'$ of $G$ and a continuous homomorphism $S:W'\to H$ such that, for $x,y\in (W')^\omega$, if $\lim_px(p)y(p)^{-1}=1_G$, then we have
$$x E_*(G)y\iff S^\omega(x)E_*(H)S^\omega(y).$$

Let $V\ni 1_G$ be an open subset of $G$ witnessing that $G$ is uniformly NSS. For each $k\in\omega$, define $V_{k}=\{g\in W': d_G(1_G,g)<2^{-k}\}$.

We claim that there exists a $k_0\in\omega$ such that
$$\forall (g_q)\in V_{k_0}^\omega\,(\lim_q S(g_q)=1_H\Rightarrow\lim_q g_q=1_G).$$

If not, then for each $k\in\omega$, there are a $j(k)\in\omega $ and a sequence $(g_{k,q})$ of elements of $V_k$ such that $\lim_q S(g_{k,q})=1_H$ and $g_{k,q}\notin V_{j(k)}$
for all $q\in\omega$. Then by the definition of $V$, there exists an $n_k>0$ such that
$$\forall g\in G\,(g\notin V_{j(k)}\Rightarrow\exists m\le n_k\,(g^m\notin V)).$$
Pick a large enough $q_k\in\omega$ so that $d_H(1_H,S(g_{k,q_k}))<2^{-k}/n_k$. Then we can find an $m_k\le n_k$ with ${g^{m_k}_{k,q_k}}\notin V$.

Put $M_{-1}=0$ and $M_k=m_0+\cdots+m_k$ for $k\in\omega$. Let $x\in G^\omega$ be defined as $x(p)=g_{k,q_k}$ for $p=M_{k-1}+i$ with $0\leq i<m_k$.
It is clear that $\lim_px(p)=1_G$. This implies that
$$x E_*(G)1_{G^\omega}\iff S^\omega(x)E_*(H)S^\omega(1_{G^\omega}).$$

For any $k\in\omega$, we have
$$x(M_{k-1})x(M_{k-1}+1)\cdots x(M_k-1)=g_{k,q_k}^{m_k}\notin V.$$
By Lemma~\ref{Lemm basic Tsi property}(2), we see that $x E_*(G)1_{G^\omega}$ fails.

For any $p<p'$, if $M_{k-1}\leq p< M_k$ and $M_{k'-1}\leq p'<M_{k'}$, then we have
$$\begin{aligned}
&d_H(1_{H^\omega}, S^\omega(x))|_{[p,p']}\leq\sum_{i=p}^{p'}d_H(1_H,S(x(i)))\cr
=&\sum_{i=p}^{M_k-1}d_H(1_H,S(g_{k,q_k}))+\sum_{i=k+1}^{k'-1}m_{i}d_H(1_H,S(g_{i,q_{i}}))+\sum_{i=M_{k'-1}}^{p'}d_H(1_H,S(g_{k',q_{k'}}))\cr
<&2^{-k}m_k/n_k+\sum_{i=k+1}^{k'-1}2^{-i}m_i/n_i+2^{-k'}m_{k'}/n_{k'}\cr
<&3\cdot 2^{-k}.
  \end{aligned}$$
This implies that $\lim_p\sup_{p\leq p'}d_H(1_H, S^\omega(x))|_{[p,p']}=0$. Then it follows from Lemma~\ref{Lemm basic Tsi property}(2) that $S^\omega(x) E_*(H)S^\omega(1_{G^\omega})$. A contradiction!

Now put $U=V_{k_0+2}$. Then $\overline{U}\subseteq V_{k_0+1}=V_{k_0+1}^{-1}$, so $\overline{U}^{-1}\overline{U}\subseteq V_{k_0}$. For $g\in\overline{U}$ and $(g_p)\in\overline{U}^\omega$, since $g_p^{-1}g\in V_{k_0}$ for each $p$, we have
$$\lim_pS(g_p)=S(g)\Rightarrow\lim_pS(g_p^{-1}g)=1_H\Rightarrow\lim_pg_p^{-1}g=1_G\Rightarrow\lim_pg_p=g.$$
This implies that $S\upharpoonright\overline{U}$ is a topological embedding.

Then we only need to show that $S(\overline{U})$ is closed. Let $(g_p)\in\overline{U}^\omega$ and $h\in H$ with $\lim_pS(g_p)=h$. Assume for contradiction that $(g_p)$ is not converges in $\overline{U}$. Then it is not $d_G$-Cauchy, we can find two strictly increasing sequences of natural numbers $(p_k),(q_k)$ with $p_k<q_k$ such that $\lim_kd_G(g_{p_k},g_{q_k})\ne 0$. By $\lim_kS(g_{p_k}^{-1}g_{q_k})=\lim_kS(g_{p_k})^{-1}S(g_{q_k})=1_H$, we have $\lim_kg_{p_k}^{-1}g_{q_k}=1_G$. We get a contradiction as desired. Therefore, $\lim_pg_p=g$ for some $g\in\overline{U}$, and hence $h=\lim_pS(g_p)=S(g)\in S(\overline{U})$. So $S(\overline{U})$ is closed.

Finally, if $G=W$ and the interval $[0,1]$ can be embedded into $H$, by lemmas~\ref{locally homeomorphism} and~\ref{converse}, we see that the converse is also true.
\end{proof}

\begin{corollary}\label{not uninfor NSS}
Let $G,H$ be two TSI Polish groups. Assume that $G$ is uniformly NSS and $H$ is strongly NSS. If $E(G)\leq_B E(H)$, then $G$ is strongly NSS too.
\end{corollary}

\begin{proof}
By Theorem~\ref{unform nss}, we can get an open subset $U\ni1_G$ of $G$ and a homeomorphism $S:\overline{U}\to H$ such that $S(gg')=S(g)S(g')$ for $g,g',gg'\in \overline{ U}$. Let $V$ be an unenclosed set of $H$. It is clear that $S^{-1}(V)\cap U$ is an unenclosed set of $G$, so $G$ is strongly NSS.
\end{proof}

\begin{remark}
The assumption of uniformly NSS in Theorem~\ref{unform nss} and Corollary~\ref{not uninfor NSS} can not be avoided. For example, we have $E(\T_p)\le_BE(\T)$, where $\T$ is strongly NSS but $\T_p$ is not (see~\cite[Theorem 6.24]{DZ}).
\end{remark}

In the rest of this article, we attempt to study some examples that are induced by totally disconnected TSI Polish groups.

\begin{example}
Let $G=\{v\in c_0:\forall n\,(v(n)\in 2^{-n}\mathbb{Z})\}$. Then $G$ is uniformly NSS but not strongly NSS with the subspace topology inherited from $c_0$. By Theorem~\ref{pro strong nss banach}, $l_p$ is strongly NSS. So $E(G)\nleq_B E(l_p)$ for $p\in[1,+\infty)$.
\end{example}

Recall that a sequence $(v_j)$ in a Banach space $X$ is called a \textit{Schauder basis} of $X$ if
$$\forall v\in X\,\exists!(t_j)\in\R^\omega\,(v=\sum_jt_jv_j).$$
Two sequences $(v_j)$ and $(v'_j)$ are called \textit{equivalent} in $X$ if
$$\sum_{j}t_jv_j\mbox{ converges}\iff\sum_{j}t_jv'_j\mbox{ converges}.$$

For $p\in[1,+\infty)$ and $a\in c_0$, we define
$$I_a=\{n\in\omega: a(n)\neq 0\},\quad A_{p,a}=\{v\in l_p:\forall n\,(v(n)\in a(n)\mathbb{Z})\}.$$
If $I_a$ is infinite, then $A_{p,a}$ equipped the relative topology inherited from $l_p$ is a totally disconnected, strongly NSS abelian Polish group, but is not non-archimedean.
In particular, we put $d(n)=2^{-n}$, and let
$$A_p=A_{p,d}=\{v\in l_p:\forall n\,(v(n)\in 2^{-n}\mathbb Z)\}.$$

Let $e_j=(0,\dots,0,\overset{j}1,0,\dots)$, then $(e_j)$ is a Schauder basis of $l_p$.

\begin{theorem}
For $p,q\in[1,+\infty)$ and $a\in c_0$, the following hold:
\begin{enumerate}
  \item [(1)] if $I_a$ is a nonempty finite set, then $E(A_{p,a})\sim_B E_0$;
  \item [(2)] if $I_a$ is infinite, then $E(A_{p,a})\sim_B E(A_p)$;
  \item [(3)] $E(A_p)<_B E(l_q)$;
  \item [(4)] $E(A_p)\leq_B E(l_q)\iff p=q$; and
  \item [(5)] $E(A_p)\leq_B E(A_q)\iff p=q$.
\end{enumerate}
\end{theorem}

\begin{proof}
  (1) If $I_a$ is a nonempty finite set, then $A_{p,a}$ is a nontrivial countable discrete group. So $E(A_{p,a})\sim_B E_0$ (cf.~\cite[Theorem 3.5(1)]{DZ}).

(2) Suppose that $I_a$ is infinite. We can find a strictly increasing sequence of $j(n)\in\omega$ for each $n\in I_a$ such that $2^{-j(n)}\le|a(n)|$.
Clearly, there is an $m_n\in\omega$ with
$$|a(n)|\leq m_n2^{-j(n)}\leq 2|a(n)|.$$

Let $\varphi:A_{p,a}\rightarrow A_p$ be defined as
$$\varphi(v)(k)=\left\{\begin{array}{ll}
\frac{m_n v(n)}{2^{j(n)}a(n)}, & k=j(n), n\in I_a,\cr
0, & \mbox{otherwise}.\cr
  \end{array}\right.$$
It is obvious that $\varphi$ is a continuous homomorphism.
For any $v\in A_{p,a}$, we have
$$\Vert v\Vert_p=\left(\sum_{n\in I_a} |v(n)|^p\right)^{\frac{1}{p}}
\leq\left(\sum_{n\in I_a}\left|\frac{m_n v(n)}{2^{j(n)}a(n)}\right|^p\right)^{\frac{1}{p}}=\Vert\varphi(v)\Vert_p\leq 2 \Vert v\Vert_p.$$
Thus $\varphi$ is a topological group isomorphism from $A_{p,a}$ onto $\varphi(A_{p,a})$. It follows from Proposition~\ref{pro closed group} that $E(A_{p,a})\leq_B E(A_p)$. Similarly, $E(A_p)\leq_B E(A_{p,a})$.

(3) Since $A_p$ is a closed subgroup of $l_p$, we get that $E(A_p)\leq_B E(l_p)$. By Corollary~\ref{cor totall disc}, we have that $E(A_p)<_B E(l_p)$.

(4) The $(\Leftarrow)$ part is from (3). We prove the $(\Rightarrow)$ part as follows:

Suppose that $E(A_p)\leq_B E(l_q)$. Then by Theorem~\ref{unform nss}, there exist an open subgroup $W$ of $A_p$ and a continuous homomorphism $S:W\to l_q$
such that $S\upharpoonright\overline{U}:\overline{U}\rightarrow S(\overline{U})$ is a homeomorphism and $S(\overline{U})$ is closed in $l_q$, where $U\subseteq W$ is an open neighborhood of $0$.

For $k\in\omega$, let $V_k=\{v\in A_p:\Vert v\Vert_p\leq 2^{-k}\}$. Since $S$ is continuous, we can find a integer $k_0>0$ so that
$$V_{k_0}\subseteq U,\quad\forall v\in V_{k_0}\,(\Vert S(v)\Vert_q\leq 1).$$
For any $k\in\omega$, we define $v_k\in V_{k_0}$ as
$$v_k(n)=\left\{\begin{array}{ll} 2^{-k_0-1},&n=k+k_0,\cr 0, & \mbox{otherwise.}\end{array}\right.$$
Note that for any $k\neq k'$, we have
$$\Vert v_k\Vert_p=2^{-k_0-1}<\Vert v_k-v_{k'}\Vert_p=2^{1/p}2^{-k_0-1}\leq 2^{-k_0},$$
and thus both $v_k$ and $v_k-v_{k'}$ are in $V_{k_0}$. So
$$\sup_{k}\Vert S(v_k)\Vert_q\leq 1,\quad\sup_{k\neq k'}\Vert S(v_k)-S(v_{k'})\Vert_q\leq 1.$$
Since $S\upharpoonright\overline{U}$ is a homeomorphism, there are $D_1,D_2>0$ such that
$$\inf_k\Vert S(v_k)\Vert_q>D_1,\quad\inf_{k\neq k'}\Vert S(v_k)-S(v_{k'})\Vert_q>D_2. $$

Since $\sup_{k}\Vert S(v_k)\Vert_q\leq 1$, by the compactness of $[-1,1]^\omega$, there are a subsequence $(S(v_{k_i}))$ of $(S(v_k))$ and a $w\in[-1,1]^\omega$ such that
$$\forall n\,(\lim_iS(v_{k_i})(n)=w(n)).$$
It is clear that $\|w\|_q\le 1$, and thus $w\in l_q$. Put $w_i=S(v_{k_{2i}})-S(v_{k_{2i+1}})\in l_q$ for each $i\in\omega$. We have
$$D_2<\inf_i\Vert w_i\Vert_q\leq\sup_i\Vert w_i\Vert_q\leq 1,\quad \forall n\,(\lim_i w_i(n)=0).$$
It follows from~\cite[Proposition l.a.12]{LT} that there is a subsequence $(w_{i_j})$ of $(w_i)$ which is equivalent to a block basis (see~\cite[Definition l.a.l0]{LT}) $(w'_j)$ of $(e_j)$. In other words, we get that, for any $(t_j)\in\R^\omega$,
$$\sum_{j}t_jw_{i_j}\mbox{ converges}\iff\sum_{j}t_jw'_j\mbox{ converges}.$$
Then it is routine to check that $0<\inf_j\Vert w'_j\Vert_q\leq\sup_j\Vert w'_j\Vert_q<\infty.$
By similar arguments in proof of~\cite[Proposition  2.a.l(i)]{LT}, we see that $(w'_j)$ is equivalent to $(e_j)$ in $l_q$. Therefore, $(w_{i_j})$ is equivalent to $(e_j)$ in $l_q$.

Let $x(j)=v_{k_{2i_j}}-v_{k_{2i_j+1}}$ for each $j\in\omega$. Then $x\in V_{k_0}^\omega\subseteq U^\omega\subseteq W^\omega$.
It is trivial to check that $(x(j))$ is equivalent to $(e_j)$ in $l_p$.
Note that $2^{-k_0-1}<\Vert x(j)\Vert_p\leq 2^{-k_0}$ and $D_2\leq \Vert w_{i_j}\Vert_q\leq 1$ for $j\in\omega$.
Thus for any $(t_j)\in\R^\omega$, we have
$$\lim_jt_j{x}(j)=0\iff\lim_jt_j=0\iff \lim_jt_jw_{i_j}=0.$$
Put $X =\{(t_j)\in\R^\omega:\forall j\,(t_j\in 2^{-j}\Z)\}$. Clearly, for any $(t_j)\in X$, $\sum_jt_jx(j)$ is in $A_p$ if it converges.
By $x(j)\in V_{k_0}\subseteq U\subseteq W$, we have $S(t_j{x}(j))=t_jw_{i_j}$ for each $(t_j)\in X$.
From Lemma~\ref{locally homeomorphism} and Proposition~\ref{equiv for sharp}, we can get the following conclusion that
$$\sum_{j}t_j{x}(j)\mbox{ converges}\iff \sum_{j}t_jw_{i_j}\mbox{ converges},$$
because both sides of the above formula imply $\lim_jt_jx(j)=0$. So
$$\forall (t_j)\in X\,((t_j)\in l_p\iff (t_j)\in l_q),$$
this gives that $p=q$.

(5) It follows immediately from (3) and (4).
\end{proof}

\end{document}